\documentclass{article}
\usepackage{graphicx} 
\usepackage{amsmath} 
\usepackage{amssymb} 
\usepackage{color}
\usepackage{mathrsfs}
\usepackage{verbatim}
\usepackage{xypic}

\newenvironment{graffa}{\left\{ \begin{array}{l}}{\end{array} \right.}
\newenvironment{graffa3}{\left\{ \begin{array}{lcc}}{\end{array} \right.}
\newenvironment{dimostr}{\textbf{Proof:}}{\hfill\rule{1ex}{1ex}}
\newtheorem{teo}{Theorem}

\newcommand{\Id}{\mathrm{Id}}
\newcommand{\nei}{neighborhood}

\newcommand{\Kil}{\mathrm{Kil}}
\newcommand{\g}{{\mathbf g}}

\renewcommand{\b}[1]{{\bf #1}}
\newcommand{\funz}[5]{#1 : \begin{tabular}{ccl}
 #2 &$\rightarrow$& #3 \\
 #4 & $\mapsto$& #5   \end{tabular}}

\newcommand{\mfunz}[5]{\funz{$#1$}{$#2$}{$#3$}{$#4$}{$#5$}}
\newcommand{\mmfunz}[5]{\begin{center}
\funz{$\displaystyle{#1}$}{$\displaystyle{#2}$}{$\displaystyle{#3}$}{$\displaystyle{#4}$}{$\displaystyle{#5}$}
\end{center}}
\newcommand{\fz}[3]{#1:\, #2 \rightarrow #3}

\newcommand{\immagine}[3][17]{
\begin{figure}[hbpt]
\begin{center}
\includegraphics[width=#1cm]{#2.eps}
\caption{#3}
\label{fig:#2}
\end{center}
\end{figure}}

\newcommand{\immintro}[3][15]{
\begin{figure}[hbpt]
\begin{center}
\includegraphics[width=#1cm]{#2.eps}
\caption{#3}
\label{intro:#2}
\end{center}
\end{figure}}

\renewcommand{\span}[1]{\mathrm{span}\Pg{#1}}
\renewcommand{\cosh}{\mathrm{Cosh}}
\renewcommand{\sinh}{\mathrm{Sinh}}
\renewcommand{\tanh}{\mathrm{Tanh}}

\newcommand{\EXP}{\mathrm{Exp}}
\newcommand{\Exp}{\mathrm{Exp}}

\newcommand{\Tr}{\mathrm{Tr}}
\newcommand{\Mat}{\mathrm{Mat}}
\newcommand{\Pt}[1]{\left( #1 \right)}
\newcommand{\Pg}[1]{\left\{ #1 \right\}}
\newcommand{\Pq}[1]{\left[ #1 \right] }
\newcommand{\hp}{hypothesis}
\newcommand{\hps}{hypotheses}

\renewcommand{\Re}[1]{\mathrm{Re}\Pt{#1}}
\renewcommand{\Im}[1]{\mathrm{Im}\Pt{#1}}

\newcommand{\bbibitem}{\bibitem}

\newcommand{\llabel}[1]{{\label{#1}}}

\newcommand{\ffoot}[1]{}

\renewcommand{\r}[1]{(\ref{#1})}
\newcommand{\bi}{\begin{itemize}}
\newcommand{\ei}{\end{itemize}}
\newcommand{\bd}{\begin{description}}
\newcommand{\ed}{\end{description}}
\newcommand{\be}{\begin{enumerate}}
\newcommand{\ee}{\end{enumerate}}

\renewcommand{\i}{\item}

\newcommand{\bqn}{\begin{eqnarray}}
\newcommand{\eqn}{\end{eqnarray}}
\newcommand{\eqnn}{\nonumber\end{eqnarray}}
\newcommand{\eqnl}[1]{\llabel{#1}\end{eqnarray}}

\newcommand{\nn}{\nonumber}

\newcommand{\ba}[1]{\begin{array}{#1}}
\newcommand{\ea}{\end{array}}

\newcommand{\R}{\mathbb{R}}
\newcommand{\C}{\mathbb{C}}

\newcommand{\fine}{\end{document}}

\def \trait (#1) (#2) (#3){\vrule width #1pt height #2pt depth #3pt}
\def \qed{\hfill
        \trait (0.1) (6) (0)
        \trait (6) (0.1) (0)
        \kern-6pt   
        \trait (6) (6) (-5.9)
        \trait (0.1) (6) (0)
\medskip}

\newtheorem{ml}{\bf Lemma}
\newtheorem{Theorem}{\bf Theorem}

\newtheorem{mrem}{\bf \underline{{\sl Remark}}}
\newtheorem{mcc}{\bf Corollary}
\newtheorem{Definition}{\bf Definition}
\newtheorem{mpr}{\bf Proposition}
\newtheorem{mproperty}{\bf Property}

\newcommand{\bt}{\begin{Theorem}}
\newcommand{\et}{\end{Theorem}}
\newcommand{\bl}{\begin{ml}}
\newcommand{\el}{\end{ml}}
\newcommand{\bp}{\begin{mpr}}
\newcommand{\ep}{\end{mpr}}
\newcommand{\bc}{\begin{mcc}}

\newcommand{\bproperty}{\begin{mproperty}}
\newcommand{\eproperty}{\end{mproperty}}

\newcommand{\ec}{\end{mcc}}
\newcommand{\bdeff}{\begin{Definition}}
\newcommand{\edeff}{\end{Definition}}
\newcommand{\brem}{\begin{mrem}\rm}
\newcommand{\erem}{\end{mrem}}

\newcommand{\proof}{{\bf Proof. }}
\newcommand{\Lam}{\Lambda}
\newcommand{\lam}{\lambda}
\newcommand{\al}{\alpha}
\newcommand{\de}{\delta}
\newcommand{\om}{\omega}

\renewcommand{\th}{\theta}

\newcommand{\Z}{{\mathbb Z}}

\newcommand{\D}{{\cal D}}
\renewcommand{\l}{\mbox{{\footnotesize \bf L}}}

\newcommand{\e}{\mbox{e}}

\newcommand{\eproof}{\hfill $\blacksquare$}

\renewcommand{\k}{{\mbox{\bf k}}}
\newcommand{\p}{{\mbox{\bf p}}}

\newcommand{\ga}{\gamma}

\begin{document} 

\begin{center} \noindent
{\LARGE{\sl{\bf Invariant Carnot-Caratheodory metrics on $S^3$, $SO(3)$, 
$SL(2)$ and lens spaces}}}
%\vskip 1cm
%\today
\vskip 1cm
Ugo Boscain

{\footnotesize LE2i,  CNRS  UMR5158, 
Universit\'e de Bourgogne,
9, avenue Alain Savary - BP 47870,   
21078 Dijon CEDEX, France}

and

{\footnotesize SISSA, via Beirut 2-4 34014 Trieste, Italy - {\tt boscain@sissa.it}}\\
\vspace*{.5cm}

Francesco Rossi\\
{\footnotesize SISSA, via Beirut 2-4 34014 Trieste, Italy - {\tt rossifr@sissa.it}}

\end{center}

\vspace{.5cm} \noindent \rm

\begin{quotation}
\noindent  {\bf Abstract}
In this paper we study the Carnot-Caratheodory metrics on $SU(2)\simeq S^3$, $SO(3)$ and $SL(2)$ induced by their Cartan decomposition and by the Killing form. Besides computing explicitly geodesics and conjugate loci, we compute the cut loci (globally) and we give the expression of the Carnot-Caratheodory distance as the inverse of an elementary function. We then prove that the metric given on $SU(2)$ projects on the so called lens spaces $L(p,q)$. Also for lens spaces, we compute the cut loci (globally). 

For $SU(2)$ the cut locus is a maximal circle without one point. In all other cases the cut locus is a stratified set.
To our knowledge, this is the first explicit computation of the whole cut locus in sub-Riemannian geometry, except for the trivial case of the Heisenberg group.

\ffoot{Our results have connections with the construction of the heat kernels for the hypoelliptic operators associated to the sub-Riemannian structure.}
\end{quotation}

%%%%%%%%%%%%%%%%%%%%%%%%%
\vskip 1cm
\begin{center}
PREPRINT SISSA 58/2007/M
\end{center}
\vskip 1cm

\newpage
%%%%%%%%%%%%%%%%%%%%%%%%%
\section{Introduction}
In this paper we study the global structure of the cut locus (set of points reached optimally by more than one geodesic) for the simplest sub-Riemannian structures on three dimensional simple  Lie groups (i.e. $SU(2)$, $SO(3)$, and  $SL(2)$) namely, the left-invariant sub-Riemannian structure induced by their Cartan decomposition and by the Killing form. 

Let $G$ be a simple real Lie group of matrices with associated Lie algebra $\l$ and Killing form $\Kil(\cdot,\cdot)$. Let $\l=\k \oplus \p$ be its Cartan decomposition with the usual commutation relations 
$[\k,\k]\subseteq\k,~~ [\p,\p]\subseteq\k,~~[\k,\p]\subseteq\p$. If $\l$ is non compact we also require $\k$ to be the maximal compact subalgebra of $\l$.
The most natural left-invariant sub-Riemannian structure that one can define on $G$ is the one in which the distribution is generated by left translations of $\p$  and the sub-Riemannian metric  $<\cdot,\cdot>$ at the identity is generated by a scalar multiple of the Killing form restricted to $\p$. The scalar must be chosen positive or negative in such a way that the scalar product is positive definite. We call $G$, endowed with such a sub-Riemannian structure, a {\bf $\k\oplus\p$ sub-Riemannian manifold}. 

$\k\oplus\p$ sub-Riemannian manifolds  have very special features: there are no strict abnormal minimizers and the Hamiltonian system given by the Pontryagin Maximum Principle is integrable in terms of elementary functions (products of exponentials). More precisely, if we write the distribution at a point $g\in G$ as $\Delta(g)=g\p$, we have the following expression for geodesics parametrized by arclength, starting at time zero from $g_0$ (\cite{agra-book,q2,brokko-cdc99,jurd-MCT,jurd-SCL}):
\bqn
g(t)=g_0e^{(A_k+A_p)t}e^{-A_k t},
\label{uno}
\eqn
where $A_k\in\k$, $A_p\in \p$, and we have $<A_p,A_p>=1$.
Thanks to left-invariance,  with no loss of generality we can always assume $g_0$ to be the identity and we will do so all along the paper.

In all three-dimensional cases (i.e. $SU(2)$, $SO(3)$ and $SL(2)$), $\p$ has 
dimension 2, while $\k$ has dimension 1.
Writing $\p=\span{p_1,p_2}$ where $\{p_1,p_2\}$ is an orthonormal frame for the sub-Riemannian 
structure (i.e. $<p_i,p_j>=\delta_{ij}$) and $\k=\span{k}$, we can write $A_p=\cos (\th) 
p_1+\sin(\th) p_2$ and $A_k=c k$ with $\th\in \R/2\pi,~ c\in\R$. The map associating to the triple $(\th,c,t)$ the final point of the corresponding geodesic starting from the identity, is called the {\it exponential map}:
\mmfunz{\Exp}{S^1\times\R\times\R^+}{G}{(\th,c,t)}{\Exp(\th,c,t)=e^{(A_k+A_p)t}e^{-A_k t}.}

For three dimensional $\k\oplus\p$ sub-Riemannian manifolds, the local structure of the sub-Riemannian spheres, cut loci and conjugate loci starting from the identity has been described by Agrachev (unpublished) and, due to cylindrical symmetry of the Killing form in the $\p$ subspace, it is very similar to the one of the Heisenberg group. Indeed, locally, the cut locus coincides with the first conjugate locus (i.e. the set where local optimality is lost) and it is made by two connected one-dimensional manifolds adjacent to the identity and transversal to the distribution, see Figure \ref{intro:Heisenberg}.
\immintro[11]{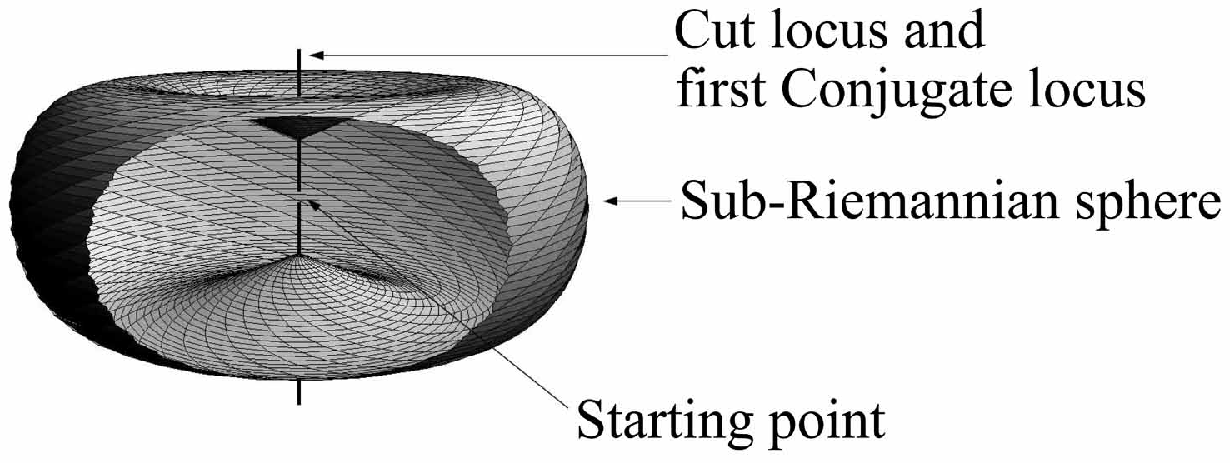}{Local structure of sub-Riemannian spheres, cut and conjugate loci for 3-dim $\k\oplus\p$ sub-Riemannian manifolds.}

However the global structure of the cut locus was still unknown. Indeed, to our knowledge, no global structure of the cut locus is known in sub-Riemannian geometry apart from the one of the Heisenberg group.

The main result of our paper is the following:
\bt
Let $K_\Id$ be the cut locus starting from the identity. We have the following:
\bi
\i for $SU(2)$, $K_\Id$ is a maximal circle $S^1$ without one point (the identity).  
\i for $SO(3)$, $K_\Id$ is a stratified set made by two manifolds glued in one point. The first manifold is $\R \mathbb{P}^2$, the second manifold is a maximal circle $S^1$ without one point (the identity).
\i for $SL(2)$,  $K_\Id$ is a stratified set made by two manifolds glued in one point. The first manifold is $\R^2$, the second manifold is a circle $S^1$ without one point (the identity).  
\ei
\et

For all cases the one dimensional strata contains the cut locus appearing in the local analysis.
\immintro{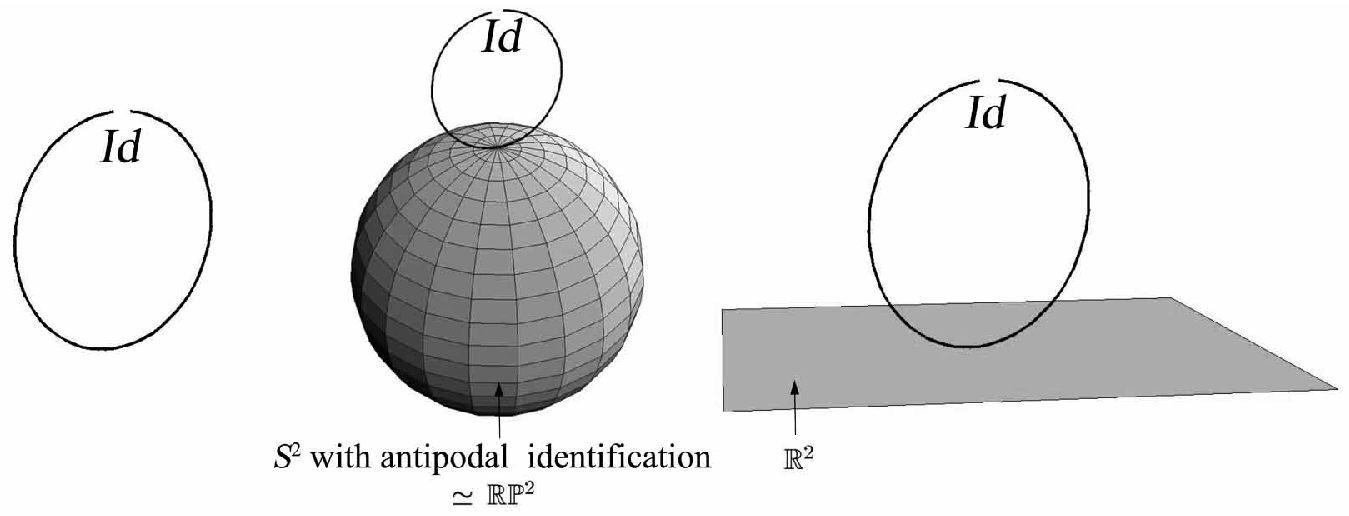}{The cut loci for the $\k\oplus\p$ sub-Riemannian manifolds $SU(2)$, $SO(3)$, $SL(2)$.}

Notice that $\k\oplus\p$ sub-Riemannian manifold $SU(2)$ has the structure of CR manifold and it is a tight structure \cite{boggess,elias}.

Once the cut locus is computed, one can provide the expression of the sub-Riemannian distance from the identity. The following theorem gives the sub-Riemannian distance for $SU(2)$. The proof, given in Section \ref{ss-SU2dist}, can be adapted to get similar results in the cases of $SO(3)$ and $SL(2)$.
\bt
\llabel{t-distance}
Let $g=\Pt{\ba{cc}\al& \beta\\-\overline{\beta}& \overline{\al}\ea}\in SU(2)$, i.e. $\al,\beta\in\C,~|\al|^2+|\beta|^2=1$ . Its sub-Riemannian distance from $\Id$ is
\bqn
\llabel{eq-SU2dist}
d(g,\Id)=\begin{graffa3}
2\sqrt{\arg(\al)\Pt{2\pi-\arg(\al)}}& \mbox{~if~} & \beta=0\\
\psi(\al) &\mbox{~if~} & \beta\neq 0
\end{graffa3},
\eqn
where $\arg\Pt{\al}\in[0,2\pi]$ and $\psi(\al)=t$ the unique solution of
\bqn
\llabel{eq-SU2distsystem}
\begin{graffa}
-\frac{ct}{2}+\arctan\Pt{\frac{c}{\sqrt{1+c^2}}\tan\Pt{\frac{\sqrt{1+c^2}t}{2}}}=\arg(\al)\\
\frac{\sin\Pt{\frac{\sqrt{1+c^2}t}{2}}}{\sqrt{1+c^2}}=\sqrt{1-|\al|^2}\\
t\in\Pt{0,\frac{2\pi}{\sqrt{1+c^2}}}
\end{graffa}.
\eqn
\et

This theorem and its analogs for $SO(3)$ and $SL(2)$ are useful to give estimates for the fundamental solutions of the hypoelliptic heat equation induced by the sub-Riemannian structure (\cite{beals,polidoro,gaveau,hormander}). Moreover this Theorem can be seen as the answer, in the case of $SU(2)$, to the question (formulated in \cite{brokko-cdc99}) about the possibility of inverting the matrix equation \r{uno}, i.e., for every matrix $g\in SU(2)$, find a matrix  $A=A_k+A_p$, with $<A_p,A_p>=1$, solution to   the equation $g=g_0e^{(A_k+A_p)t}e^{-A_k t}$. If $\beta\neq 0$ then this equation has one and only one solution, otherwise it has more than one solution (indeed infinitely many, see Sections  \ref{s-gruppi} and \ref{s-cut}).

Then we study the most natural sub-Riemannian structures on the lens spaces $L(p,q)$ induced by the one on $SU(2)$.
The lens space $L(p,q)$ (with $p,q$ coprime integers, $p,q\neq 0$) is the quotient of $SU(2)$ by the equivalence relation
$$\Pt{\ba{cc}\al_1& \beta_1\\-\overline{\beta_1}& \overline{\al_1}\ea}\sim\Pt{\ba{cc}\al_2& \beta_2\\-\overline{\beta_2}& \overline{\al_2}\ea}\mbox{~if~$\exists~\om\in\C$ $p$-th root of unity such that~}\Pt{\ba{c}\al_2\\\beta_2\\\ea}=
\Pt{\begin{array}{cc}
\om & 0\\
0 & \om^q
\end{array}}\Pt{\ba{c}\al_1\\\beta_1\\\ea}.$$
The lens spaces are three-dimensional manifolds, but they are neither Lie groups nor homogeneous spaces of $SU(2)$, except for the case $L(2,1)\simeq SO(3)$ \ffoot{Chiedere a Luisa}.

In the case of lens spaces we get that the cut locus is much more complicated with respect to those on $SU(2)$ and $SL(2)$. It is still a stratified set, but in general with more strata. The precise description is given in Section \ref{ss-lpcut}.

Sub-Riemannian structures on the lens space $L(4,1)$ are particularly interesting for mechanical applications and for problems of geometry of vision on the two-dimensional sphere. Indeed, $L(4,1)\simeq PTS^2$, the bundle of directions of $S^2$. These applications are the subject of a forthcoming paper.

The structure of the paper is the following: in Section \ref{s-basic} we recall the definition of sub-Riemannian manifold, we state the Pontryagin Maximum Principle (that is a first order necessary condition for optimality for problems of calculus of variations with non-holonomic constraints) and we define the cut and conjugate loci. Then we define $\k\oplus\p$ sub-Riemannian manifolds. In Section \ref{s-gruppi} we define $\k\oplus\p$ sub-Riemannian structures on $SU(2),SO(3), SL(2)$ and compute the corresponding geodesics and conjugate loci. In Section \ref{s-Lpq} we give sub-Riemannian structures on lens spaces as quotients of the $\k\oplus\p$ sub-Riemannian structure on $SU(2)$. The core of the paper is Section \ref{s-cut}, where we compute the cut loci and the sub-Riemannian distance. The general idea is the following: we first identify the prolongation of the cut locus arising locally, then we compute the part of the cut locus due to the symmetries of the problem and finally we show that there is no other cut point.

\section{Basic Definitions}
\llabel{s-basic}
\subsection{Sub-Riemannian manifold}
A $(n,m)$-sub-Riemannian manifold is a triple $(M,\Delta,{\mathbf g})$, 
where
\bi
\i $M$ is a connected smooth manifold of dimension $n$;
\i $\Delta$ is a Lie bracket generating smooth distribution of constant rank $m<n$, i.e. $\Delta$ is a smooth map that associates to $q\in M$  a $m$-dim subspace $\Delta(q)$ of $T_qM$, and $\forall~q\in M$ we have
\bqn\llabel{Hor}\span{[f_1,[\ldots[f_{k-1},f_k]\ldots]](q)~|~f_i\in\mathrm{Vec}(M)\mbox{~and~}f_i(p)\in\Delta(p)~\forall~p\in M}=T_qM.\eqn
Here $Vec(M)$ denotes the set of smooth vector fields on $M$.
\i ${\mathbf g}_q$ is a Riemannian metric on $\Delta(q)$, that is smooth 
as function of $q$.
\ei

The Lie bracket generating condition \r{Hor} is also known as H\"ormander condition.

A Lipschitz continuous curve $\ga:[0,T]\to M$ is said to be \b{horizontal} if 
$\dot\ga(t)\in\Delta(\ga(t))$ for almost every $t\in[0,T]$.
Given an horizontal curve $\ga:[0,T]\to M$, the {\it length of $\ga$} is
\bqn
l(\ga)=\int_0^T \sqrt{ \g_{\ga(t)} (\dot \ga(t),\dot \ga(t))}~dt.
\eqn
The {\it distance} induced by the sub-Riemannian structure on $M$ is the 
function
\bqn
d(q_0,q_1)=\inf \{l(\ga)\mid \ga(0)=q_0,\ga(T)=q_1, \ga\ \mathrm{horizontal}\}.
\eqnl{e-dipoi}

The \hp\ of connectedness of M and the Lie bracket generating assumption for the distribution guarantee the finiteness and the continuity of $d(\cdot,\cdot)$ with respect to the topology of $M$ (Chow's Theorem, see for instance \cite{agra-book}).

The function $d(\cdot,\cdot)$ is called the Carnot-Charateodory distance and gives to $M$ the structure of metric space (see \cite{bellaiche,gromov}).

It is a standard fact that $l(\ga)$ is invariant under reparameterization of the curve $\ga$.
Moreover, if an admissible curve $\ga$ minimizes the so-called {\it energy functional}
$$ E(\ga)=\int_0^T {\g}_{\ga(t)}(\dot \ga(t),\dot \ga(t))~dt. $$
with $T$ fixed (and fixed initial and final point), then $v=\sqrt{\g_{\ga(t)}(\dot \ga(t),\dot \ga(t))}$
is constant and $\ga$ is also a minimizer of $l(\cdot)$.
On the other side a minimizer $\ga$ of $l(\cdot)$ such that  $v$ is constant is a minimizer of $E(\cdot)$ with $T=l(\ga)/v$.

A {\it geodesic} for  the sub-Riemannian manifold  is a curve $\ga:[0,T]\to M$ such that for every sufficiently small interval $[t_1,t_2]\subset [0,T]$, $\ga_{|_{[t_1,t_2]}}$ is a minimizer of $E(\cdot)$.
A geodesic for which $\g_{\ga(t)}(\dot \ga(t),\dot \ga(t))$  is (constantly) equal to one is said to be parameterized by arclength.

Locally, the pair $(\Delta,{\mathbf g})$ can be given by assigning a set of $m$ smooth vector fields that are orthonormal for ${\mathbf g}$, i.e.  
\bqn
\Delta(q)=Span\{F_1(q),\dots,F_m(q)\}, ~~~{\mathbf g}_q(F_i(q),F_j(q))=\delta_{ij}.
\eqnl{trivializable}
When  $(\Delta,{\mathbf g})$ can be defined as in \r{trivializable} by $m$ vector fields defined globally, we say that the sub-Riemannian manifold is {\it trivializable}. 

%%%%%%%%%%%%%%%%%%%%%%%%%%%%%%%%%%%%%%%%%%%%%%%%%%%%%%%%%%%%%%%%%%%%

Given a $(n,m)$- trivializable sub-Riemannian manifold, the problem of finding a curve minimizing the energy between two fixed points  $q_0,q_1\in M$ is
naturally formulated as the optimal control problem 
\bqn
\dot q=\sum_{i=1}^m u_i F_i(q)\,,~~~u_i\in\R\,,
~~~\int_0^T
\sum_{i=1}^m u_i^2(t)~dt\to\min,~~q(0)=q_0,~~~q(T)=q_1.
\label{sopra}
\eqn
It is a standard fact that this optimal control problem  is equivalent to the minimum time problem with controls $u_1,\ldots, u_m$ satisfying $u_1^2+\cdots+u_m^2\leq 1$.

When the manifold is analytic and the orthonormal frame can be assigned through $m$ analytic vector fields, we say that the sub-Riemannian manifold is {\it analytic}.

In this paper we are concerned with sub-Riemannian manifolds that are trivializable and analytic since they are given in terms of left-invariant vector fields on Lie groups.
%%%%%%%%%%%%%%%%%%%%%%%%%%%%%%%%%%%%%%%%%%%%%%%%%%%%%%%%%%%%%%%%%%%%%%%%%%%%%%%%%%
%%%%%%%%%%%%%%%%%%%%%%%%%%%%%%%%%%%%%%%%%%%%%%%%%%%%%%%%%%%%%%%%%%%%%%%%%%%%%%%%%%
\subsection{First order necessary conditions, Cut locus, Conjugate locus}

Consider a trivializable $(n,m)$-sub-Riemannian manifold. Solutions to the optimal control problem 
\r{sopra} are computed via the Pontryagin Maximum Principle (PMP for short, see for instance \cite{agra-book,libro,jurd-book,pontlibro}) that is a first order
necessary condition for optimality and generalizes the Weierstra\ss \
conditions of Calculus of Variations. For each optimal curve, the PMP
provides a lift to the cotangent bundle that is a solution to a suitable
pseudo--Hamiltonian system.\\ 

\bt[Pontryagin Maximum Principle for the problem \r{sopra}]
\it Let $M$ be a $n$-dimensional smooth manifold and consider the minimization problem \r{sopra}, in the class of Lipschitz continuous curves, 
 where  $F_i$, $i=1,\ldots, m$ are smooth vector fields on $M$ and the final time $T$ is fixed.  Consider  the map $\mathscr{H}:T^\ast M\times\R\times \R^m\to \R$ defined  by
\bqn
\mathscr{H}(q,\lam,p_0,u)&:=&<\lam,\sum_{i=1}^m u_i F_i(q)>+p_0 \sum_{i=1}^m u_i^2(t).
\eqnn
If the curve $q(.):[0,T]\to M$ corresponding to the control $u(.):[0,T]\to\R^m$ 
is optimal then there exist a {never vanishing}
Lipschitz continuous {covector}
$\lam(.):t\in[0,T]\mapsto \lam(t)\in
T^\ast_{q(t)}M$ and a constant $p_0\leq 0$ such that, for a.e. $t\in 
[0,T]$:
\begin{description}
\item[i)]
$\dot q(t)=\displaystyle{\frac{\partial \mathscr{H}
}{\partial \lam}(q(t),\lam(t),p_0,u(t))}$,
\item[ii)] $\dot \lam(t)=-\displaystyle{\frac{\partial \mathscr{H}
}{\partial q}(q(t),\lam(t),p_0,u(t))}$,
\item[iii)] $\frac{\partial\mathscr{H} }{\partial u}(q(t),\lam(t),p_0,u(t))=0.$
\ed
\et
\brem
\llabel{r-postPMP}
A curve $q(.):[0,T]\to M$ satisfying the PMP is said to be an {\it extremal}. In general, an extremal may correspond to more than one pair $(\lam(.),p_0)$. 
If an extremal satisfies the PMP  with $p_0\neq 0 $, then it is called a {\it normal extremals}. If it satisfies the PMP with $p_0= 0 $ it is called an {\it abnormal extremal}. An extremal can be both normal and abnormal. For normal extremals one can normalize $p_0=-1/2$.

If an extremal satisfies the PMP only with $p_0=0$, then it is called a {\it strict abnormal extremal}. If a strict abnormal extremal is optimal, then  it is called {\it a strict abnormal minimizer}. For a deep analysis of abnormal extremals in sub-Riemannian geometry, see \cite{bonnard,trelat}.
\erem

It is well known that all normal extremals are geodesics  (see for instance \cite{agra-book}). Moreover  if there are no strict abnormal minimizers then all geodesics are normal extremals for some fixed final time $T$.
This will be always the case in this paper: indeed we are concerned with sub-Riemannian manifolds of dimension 3, defined by a pair of vector fields $F_1$ and $F_2$ such that for all $q\in M$, $Span\{F_1(q),F_2(q), [F_1(q),F_2(q)]\}=T_q M$, i.e. the so called 3-D {\it contact case}, for which there are no abnormal extremals (even non strict).

%If there are no strict abnormal minimizers, then all solutions of the PMP are geodesics and all geodesics are solution of the PMP for some fixed final time $T$ (see for instance \cite{agra-book}). This will be always the case in this paper: indeed we are concerned with sub-Riemannian manifolds of dimension 3, defined by a pair of vector fields $F_1$ and $F_2$ such that for all $q\in M$, $Span\{F_1(q),F_2(q), [F_1(q),F_2(q)]\}=T_q M$, i.e. the so called 3-D {\it contact case}, for which there are no abnormal extremals (even non strict).

In this case from {\bf iii)} one gets $u_i(t)=<\lam(t),F_i(t)>$, $i=1\ldots ,m$ and the PMP becomes  much simpler:   a curve $q(.)$ is a geodesic if and only if it is the projection on $M$ of a 
solution $(\lam(t),q(t))$ for the Hamiltonian system on $T^\ast M$ corresponding to:
\bqn
H(\lam,q)=\frac12(\sum_{i=1}^m<\lam,F_i(q)>^2),~~~~
\mbox{$q\in M$, $\lam\in T^\ast_q M$}.
\eqnn
satisfying $H(\lam(0),q(0))\neq 0$.

%%%%%%%%%%%%%%%%%%%%%%%%%%%%%%%%%%%%%%%%%%%%%%%%%%%%%%%%%%%%%%%%

\brem Notice that $H$ is constant along any given solution of the Hamiltonian  system. Moreover, $H=\frac{1}{2}$ if and only if the geodesic is parameterized by arclength. In the following, for simplicity of notation, we assume that all geodesics are defined for $t\in [0,+\infty)$.
\erem
Fix $q_0\in M$. For every $\lam_0\in T_{q_0}^\ast M$ satisfying
\bqn
H(\lam_0, q_0)=1/2
\eqnl{1/2}
and every $t>0$ define the {\it exponential map} $\Exp(\lam_0,t)$ as the projection on $M$ of the solution, 
evaluated at time $t$, of the Hamiltonian system associated with $H$,  with initial condition $\lam(0)=\lam_0$ and $q(0)=q_0$.
Notice that condition \r{1/2} defines a hypercylinder $\Lambda_{q_0}\simeq S^{m-1}\times\R^{n-m}$ in  $T_{q_0}^\ast M$.
%%%%%%%%%%%%%%%%%%%%%%%%%%%%%%%%%%%%%%%%%%%%%%%%%%%%%%%%%%%%%%%%%%%%%%%%%%%%%%%%%%%%%%%%%%
\bdeff
The {\bf conjugate locus from $q_0$} is the set $C_{q_0}$ of critical values of the map \mmfunz{\EXP}{\Lam_{q_0}\times \R^+}{M}{(\lam_0,t)}{\Exp(\lam_0,t).}

For every $\bar\lam_0\in\Lam_{q_0}$, let $t(\bar\lam_0)$ be the n-th positive time, if it exists, for which the map $(\lam_0,t)\mapsto \EXP(\lam_0,t)$ is singular at $(\bar\lam_0,t(\bar\lam_0))$. The {\bf n-th conjugate locus} from $q_0$ $C^n_{q_0}$ is the set $\{\EXP(\bar\lam_0,t(\bar\lam_0))\mid t(\bar\lam_0) \mbox{ 
exists}\}$.

The {\bf cut locus} from $q_0$ is the set $K_{q_0}$ of points reached optimally by more than one geodesic, i.e., the set
\bqn
K_{q_0}=\Pg{q\in M \mid \exists~\lam_1,\lam_2\in\Lambda_{q_0},~\lam_1\neq\lam_2,~t\in\R^+\mbox{ such that~ }\ba{l}  q=\EXP(\lam_1,t)=\EXP(\lam_2,t),\mbox{~~~and~~~}\\ \EXP(\lam_1,\cdot),
\EXP(\lam_2,\cdot)\mbox{ optimal in }[0,t]\ea}
\eqnn
\edeff
%%%%%%%%%%%%%%%%%%%%%%%%%%%%%%%%%%%%%%%%%%%%
\brem
It is a standard fact that for every $\bar\lam_0$ satisfying \r{1/2}, the set $T(\bar\lam_0)=\{\bar t>0\mid 
\EXP(\lam,t)$ is singular at $(\bar\lam_0,\bar t) \}$ is a discrete set (see for instance \cite{agra-book}).
\erem
\brem
\llabel{rem-agra}
Let $(M,\Delta,\g)$ be a sub-Riemannian manifold. Fix $q_0\in M$ and assume: {\bf i)} each point of $M$ is reached by an optimal geodesic  starting from $q_0$; {\bf ii)}  there are no abnormal minimizers.
The following facts are well known (a proof in the 3-D contact case can be found in \cite{agra-exp}). 
\bi
\i the first conjugate locus $C^1_{q_0}$ is the set of points where the geodesics starting from $q_0$   lose local optimality;
\i if $q(.)$ is a geodesic starting from $q_0$ and $\bar t$ is the first positive time such that $q(\bar t)\in K_{q_0}\cup C^1_{q_0}$, then $q(.)$ loses optimality in $\bar{t}$, i.e. it is optimal in $[0,\bar{t}]$ and not optimal in $[0,t]$ for any $t>\bar{t}$;
\i if a geodesic $q(.)$ starting from $q_0$ loses optimality at $\bar t>0$, then $q(\bar t)\in K_{q_0}\cup C^1_{q_0}$;
\ei
%The following facts are well known and can be proved as in Riemannian geometry:
%\bi
%\i the first conjugate locus is the set of points where geodesics lose local optimality;
%\i let $\ga$ be a geodesic such that $\ga(0)=q_0$ and $\bar{t}$ the first positive time (if it exists) such that $\ga(\bar{t})\in K_{q_0}$ or $\bar{t}$ is a conjugate time for $\ga$; then $\ga$ loses optimality in $\bar{t}$, i.e. it is optimal in $[0,\bar{t}]$ and not optimal in $[0,t]$ for any $t>\bar{t}$.
%\ei
As a consequence, when the first conjugate locus is included in the cut locus (as in our cases, see Section \ref{s-cut}), the cut locus is the set of points where the geodesics lose optimality.
\erem
\brem
It is well known that, while in Riemannian geometry $K_{q_0}$ is never adjacent to $q_0$, in sub-Riemannian geometry this is always the case. See \cite{Agra-sub}.
\erem

%%%%%%%%%%%%%%%%%%%%%%%%%%%%%%%%%%%%%%%%%%%%%%%%%%%%%%%%%%%%%%%%%%%%%%%%%%
%%%%%%%%%%%%%%%%%%%%%%%%%%%%%%%%%%%%%%%%%%%%%%%%%%%%%%%%%%%%%%%%%%%%%%%%%%
\subsection{$\k\oplus\p$ sub-Riemannian manifolds}
\llabel{s-k+p}
%%%%%%%%%%%%%%%%%%%%%%%%%%%%%%%%%%%%%%%%%%%%%%%%%%%%%%%%%%%%%%%%%%%%%%%%%%
%%%%%%%%%%%%%%%%%%%%%%%%%%%%%%%%%%%%%%%%%%%%%%%%%%%%%%%%%%%%%%%%%%%%%%%%%%
For the sake of simplicity in the exposition, all over the paper, when we 
deal about Lie groups and Lie algebras, we always consider that they are groups
and algebras of matrices.

Let $\l$ be a simple Lie algebra and $\Kil(X,Y)=Tr(ad_X\circ ad_Y)$ its Killing 
form. Recall that the Killing form 
defines a non-degenerate pseudo scalar product on $\l$.
In the following we recall what
we mean by a Cartan decomposition of $\l$.      
\bdeff 
A Cartan decomposition of a simple Lie algebra
$\l$ is any
decomposition of the form:
\bqn
\l=\k\oplus\p, 
 \mbox{ where } [\k,\k]\subseteq\k,~~ [\p,\p]\subseteq\k,~~
[\k,\p]\subseteq\p.  
\eqnl{cartan-dec}
\llabel{d-k+p}
\edeff
%%%%%%%%%%%%%%%%%%%%%%%%%%%%%%%%%%%%%%%%%%%%%%%%%%%%%%%%%%%%%%%%%%%%%%
\bdeff 
Let $G$ be a simple Lie group with Lie algebra $\l$. Let $\l=\k\oplus \p$ be a Cartan decomposition of $\l$. In the case in which $G$ is noncompact assume that 
$\k$ is the maximal compact subalgebra of $\l$.

On $G$, consider the distribution  $\Delta(g)=g\p$ endowed with the Riemannian metric ${\mathbf g}_g(v_1,v_2)=<g^{-1}v_1,g^{-1}v_2>$ where $<~,~>:=\al~\Kil\big|_\p(~,~)$ and $\al<0$ (resp. $\al>0$) if $G$ is compact (resp. non 
compact).

In this case we say that $(G, \Delta, {\mathbf g})$ is a $\k\oplus \p$ sub-Riemannian manifold.
\llabel{d-k+p-problem}
\edeff 

The constant $\al$ is clearly not relevant. It is chosen just to obtain
good normalizations.
\brem
In the compact (resp. noncompact) case the fact that ${\mathbf g}$ is positive definite on 
$\Delta$ is guaranteed  by the requirement $\al<0$ (resp. by the requirements $\al>0$  and $\k$ 
maximal compact 
subalgebra).
\erem

Let $\{X_j\}$ be an orthonormal frame for the subspace $\p\subset\l$, with respect to the
metric defined in Definition \ref{d-k+p-problem}. 
Then the problem of finding the minimal energy between the identity and a point $g_1\in G$ in fixed time $T$ becomes the left-invariant optimal control problem
\bqn
\dot g=g\left(\sum_j u_jX_j\right),~~~~u_j\in L^\infty(0,T)\,,
\int_{0}^{T}\sum_ju_j^2(t)~dt\to\min,~~g(0)=\Id,~~~g(T)=g_1.
\eqnn
This problem admits a solution, see for instance Chapter 5 of \cite{piccoli}.

For $\k\oplus\p$ sub-Riemannian manifolds, one can prove that strict abnormal extremals are never optimal, since the {\it Goh condition} (see \cite{agra-book}) is never satisfied. Moreover, the Hamiltonian system given by the Pontryagin Maximum Principle is integrable and the explicit expression of geodesics starting from the identity and parameterized by arclength is
\bqn
g(t)=e^{(A_k+A_p)t}e^{-A_k t},
\eqnl{k+p}
where $A_k\in\k$, $A_p\in\p$ and $<A_p,A_p>=1$. This formula is known from long time in the community. It was used  independently by Agrachev \cite{agra-ICM}, Brockett \cite{brokko-cdc99} and Kupka (oral communication). The first complete proof was written by Jurdjevic in \cite{jurd-MCT}. The proof that strict  abnormal extremals are never optimal was first written  in \cite{q2}.  See also  \cite{agra-book,jurd-SCL}.

\brem
In the 3-dimensional case, the Hamiltonian system given by the Pontryagin Maximum Principle is indeed integrable even if the cost is not built with the Killing form (bi-invariant), but  is only left-invariant. For the case of $SO(3)$ see \cite{q4}.
\erem

\section{$SU(2)$, $SO(3)$, $SL(2)$, their geodesics and their conjugate loci}
\llabel{s-gruppi}
%%%%%%%%%%%%%%%%%%%%%%%%%%%%%%%%%%%%%%%%%%%%%%%%%%%%%%%%%%%%%%%%%%%%%%%%%%%%%%%%%%%%%%%%%
In this section we fix coordinates on $SU(2)$, $SO(3)$, $SL(2)$, and we apply formula \r{k+p} in order to get the explicit expressions for geodesics and conjugate loci.

\subsection{The $\k\oplus\p$ problem on $SU(2)$}
\label{ss-SU2}

The Lie group $SU(2)$ is the group of unitary unimodular $2\times 2$ complex matrices 
$$SU(2)=\Pg{\Pt{\begin{array}{cc}
\alpha & \beta\\
-\overline{\beta} & \overline{\alpha}
\end{array}}\in\Mat(2,\C)\ |\ |\al|^2+|\beta|^2=1}.$$
It is compact and simply connected. The Lie algebra of $SU(2)$ is the algebra of antihermitian traceless $2\times 2$ complex matrices 
$$su(2)=\Pg{\Pt{\begin{array}{cc}
i\alpha & \beta\\
-\overline{\beta} & -i\alpha
\end{array}}\in\Mat(2,\C)\ |\ \alpha\in\R,\beta\in\C}.$$
A basis of $su(2)$ is $\Pg{p_1,p_2,k}$ where 
\bqn
p_1=\frac{1}{2}\Pt{\begin{array}{cc}
0 & 1\\
-1 & 0
\end{array}}\quad
p_2=\frac{1}{2}\Pt{\begin{array}{cc}
0 & i\\
i & 0
\end{array}}\quad
k=\frac{1}{2}\Pt{\begin{array}{cc}
i & 0\\
0 & -i
\end{array}},
\eqnl{pauli}
whose commutation relations are $[p_1,p_2]=k\quad[p_2,k]=p_1\quad[k,p_1]=p_2$. 
Recall that for $su(n)$ we have $\Kil(X,Y)=2n\Tr(XY)$, see \cite[p. 186, 516]{helgason}; thus for $su(2)$  $\Kil(X,Y)=4\Tr(XY)$ and, in particular, $\Kil(p_i,p_j)=-2\de_{ij}$.
The choice of the subspaces $$\k=\span{k}\qquad \p=\span{p_1,p_2}$$
provides a {\it Cartan decomposition} for $su(2)$.
Moreover, $\Pg{p_1,p_2}$ is a orthonormal frame for the inner product $<\cdot,\cdot>=-\frac{1}{2}\Kil(\cdot,\cdot)$ restricted to $\p$.

Defining $\Delta(g)=g\p$ and ${\mathbf g}_g(v_1,v_2)=< g^{-1}v_1,g^{-1}v_2>$, we have that $(SU(2),\Delta,{\mathbf g})$ is a $\k\oplus\p$ sub-Riemannian manifold.
\brem
Observe that all the $\k\oplus\p$ structures that one can define on $SU(2)$ are equivalent. For instance, one could set $\k=\span{p_1}$ and $\p=\span{p_2,k}$.
\erem

Recall that $SU(2)\simeq S^3=\Pg{\Pt{\ba{c}\al\\\beta\ea}\in\C^2\ |\ |\al|^2 + |\beta|^2=1} $ via the isomorphism
\mmfunz{\phi}{SU(2)}{S^3}{\Pt{\begin{array}{cc}
\alpha & \beta\\
-\overline{\beta} & \overline{\alpha}
\end{array}}}{\Pt{\ba{c}\al\\\beta\ea}.}
In the following we always write elements of $SU(2)$ as pairs of complex numbers.

\subsubsection{Expression of geodesics}
\llabel{ss-SU2geod}
We compute the explicit expression of geodesics using the formula \r{k+p}. Consider an initial covector $\lam=\lam(\th,c)=\cos(\th)p_1+\sin(\th)p_2+ck\in \Lambda_\Id$. The corresponding exponential map is

$$\ba{l}
\EXP(\th,c,t):=\EXP(\lam(\th,c),t)=e^{(\cos(\th)p_1+\sin(\th)p_2+ck)t} e^{-ck t}=\\\\
=\Pt{\begin{array}{c}
\frac{c\sin(\frac{ct}{2})\sin (\sqrt{1 + c^2}\frac{t}{2})}{\sqrt{1 + c^2}} + \cos(\frac{ct}{2})\cos(\sqrt{1 + c^2}\frac{t}{2})+i\Pt{\frac{c\cos(\frac{ct}{2})\sin (\sqrt{1 + c^2}\frac{t}{2})}{\sqrt{1 + c^2}} - \sin(\frac{ct}{2})\cos(\sqrt{1 + c^2}\frac{t}{2})}\\
\frac{\sin(\sqrt{1 + c^2}\frac{t}{2})}{\sqrt{1 + c^2}}
\Pt{\cos(\frac{ct}{2}+\th) + i\sin(\frac{ct}{2}+\th)}
\end{array}}.\ea$$

We have the following symmetry properties
\bi
\i \textit{cylindrical symmetry}: $$\EXP(\th,c,t)=\Pt{\begin{array}{cc}
1 & 0\\
0 & e^{i\th}
\end{array}} \EXP(0,c,t);$$
\i \textit{central symmetry}:

set $\Pt{\ba{c} \alpha\\\beta\ea}=\EXP(\th,c,t)$. We have $\EXP(\th,-c,t)=\begin{graffa3}
\Pt{\ba{cc}
\overline{\alpha}\\
e^{2i(\th-\arg(\beta))}\beta
\end{array}}&\mbox{~if~}&\beta\neq 0\\
\Pt{\ba{cc}
\overline{\alpha}\\
0
\end{array}}&\mbox{~if~}&\beta= 0\\
\end{graffa3}.$
\ei
\subsubsection{Pictures of $S^2$ and $S^3$}
\llabel{ss-SU2picture}
We recall a standard construction for representing $S^2$ in a two dimensional space and $S^3$ in a three dimensional one. For more details see e.g. \cite{weeks}. Consider $S^2\subset \R^3$ and flatten it on the equator plane, pushing the northern hemisphere down and the southern hemisphere up, getting two superimposed disks $D^2$ joined along their circular boundaries. The construction is drawn in Figure \ref{fig:S2S3}-left. Similarly, consider $S^3\subset \C^2\simeq\R^4$: it can be viewed as two superimposed balls joined along their boundaries. In this case the boundaries are two spheres $S^2$. A picture of $S^3$ is given in Figure \ref{fig:S2S3}-right.
\immagine{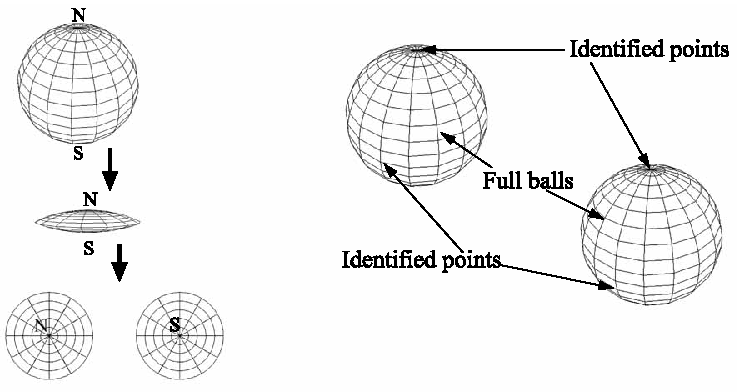}{Left: construction of the 2-dim picture of $S^2$. Right: the 3-dim picture of $S^3$.}

\subsubsection{The conjugate locus}
\label{ss-SU2conj}

Recall that all the partial derivatives of $\Exp$ evaluated in $(\th,c,t)$ lie in $T_gSU(2)=g\cdot su(2)$ with $g=\Exp(\th,c,t)$. One can easily check that the three vectors $g^{-1}\cdot \frac{\partial \EXP}{\partial \th}_{|_{(\th,c,t)}},~ g^{-1}\cdot \frac{\partial \EXP}{\partial c}_{|_{(\th,c,t)}},~g^{-1}\cdot \frac{\partial \EXP}{\partial t}_{|_{(\th,c,t)}}\in su(2)$ are linearly dependent (hence $g$ is a conjugate point) if and only if $$\sin\Pt{\sqrt{1+c^2}\frac{t}{2}}
\Pt{2\sin\Pt{\sqrt{1+c^2}\frac{t}{2}}-\sqrt{1+c^2}t\cos\Pt{\sqrt{1+c^2}\frac{t}{2}}}=0.$$

The first term is 0 if and only if $g\in e^\k=\Pg{\Pt{\begin{array}{c}
\al\\
0
\end{array}}\ |\ |\al|=1}$, while the second vanishes if and only if $\sqrt{1+c^2}\frac{t}{2}=\tan\Pt{\sqrt{1+c^2}\frac{t}{2}}$, hence we have two series of conjugate times:
\bi
\i first series: $t_{2n-1}=\frac{2n\pi}{\sqrt{1+c^2}}$ to which correspond the conjugate loci $C^{2n-1}_\Id=e^\k\setminus\Id$;
\i second series: $t_{2n}=\frac{2x_n}{\sqrt{1+c^2}}$ where $\Pg{x_1, x_2,\ldots}$ is the ordered set of the strictly positive solutions of $x=\tan(x)$, to which correspond the conjugate loci

$$ C^{2n}_\Id=\scriptstyle\Pg{\Pt{\begin{array}{c}
\frac{c\sin(x_n)}{\sqrt{1+c^2}}\Pt{\sin\Pt{\frac{cx_n}{\sqrt{1+c^2}}}+i\cos\Pt{\frac{cx_n}{\sqrt{1+c^2}}}} + \cos(x_n)\Pt{\cos\Pt{\frac{cx_n}{\sqrt{1+c^2}}}-i\sin\Pt{\frac{cx_n}{\sqrt{1+c^2}}}}\\
\\
\frac{\sin(x_n)}{\sqrt{1+c^2}} e^{i\th}
\end{array}}\mid \ba{c} c\in\R\\\th\in \R/2\pi\ea}.$$
\ei
\brem
Notice that all the geodesics have a countable number of conjugate times.
\erem

We present some images of conjugate loci (Figures \ref{fig:SU2conj1} and \ref{fig:SU2conj}). For simplicity we present an image of their section with the plane $\Re{\beta}=0$. The complete picture can be recovered using the cylindrical symmetry.
\immagine[6]{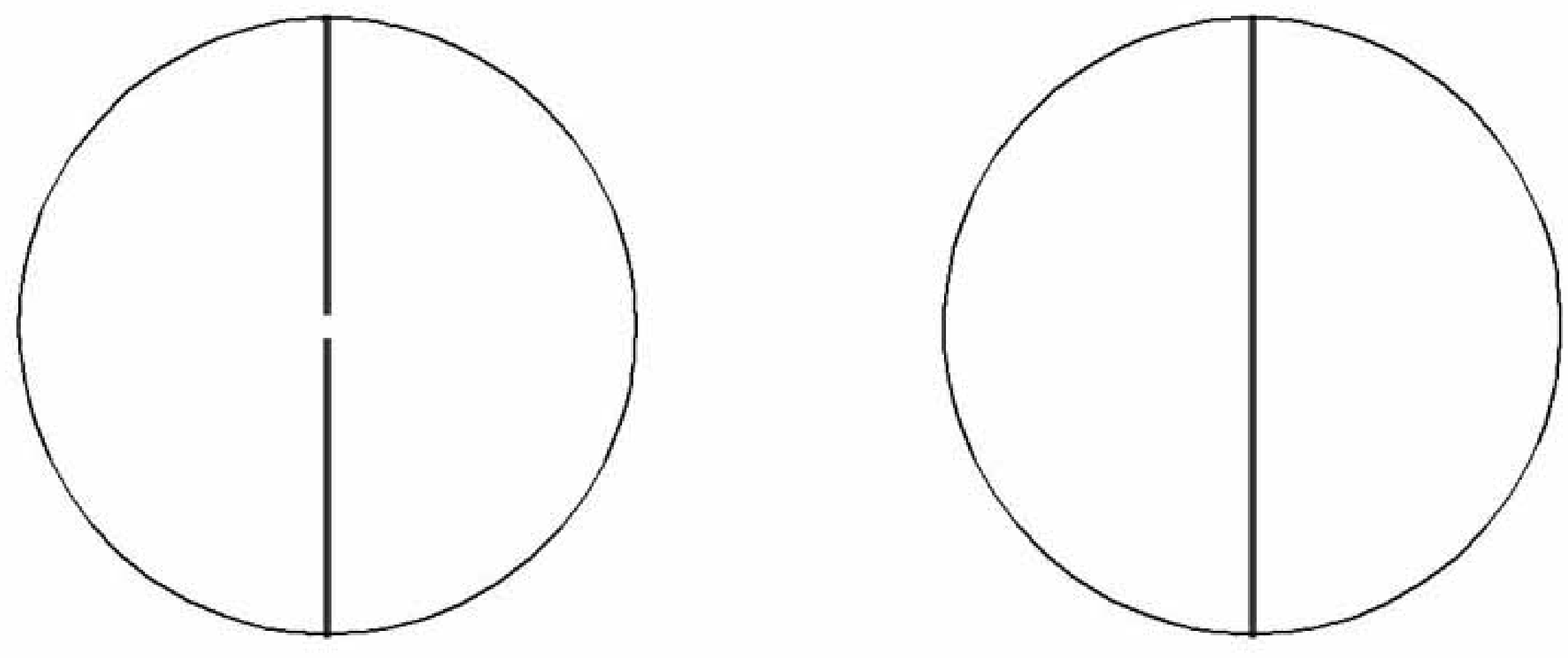}{$\k\oplus\p$ problem on $SU(2)$: projection of the odd conjugate loci.}
\immagine{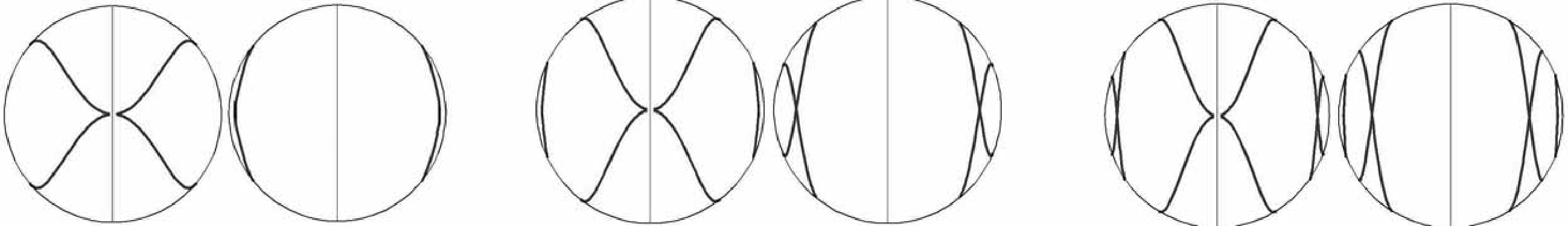}{$\k\oplus\p$ problem on $SU(2)$: projection of the $2^{nd}$, $4^{th}$ and $6^{th}$ conjugate loci.}
\brem
Notice that the second conjugate locus is a 2-dimensional submanifold of $SU(2)$, while the other even conjugate loci have self-intersections.
\erem

\subsection{The $\k\oplus\p$ problem on $SO(3)$}
\label{ss-SO3}

The Lie group $SO(3)$ is the group of special orthogonal $3\times 3$ real matrices 
$$SO(3)=\Pg{g\in\Mat(3,\R)\ |\ g g^T=\Id, \det(g)=1}.$$
It is compact and its fundamental group is $\Z_2$. The Lie algebra of $SO(3)$ is the algebra of skew-symmetric $3\times 3$ real matrices $$so(3)=\Pg{\Pt{\ba{ccc}
0 & -a & b\\
a & 0 & -c\\
-b & c & 0\ea}\in\Mat(3,\R)}.$$
A basis of $so(3)$ is $\Pg{p_1,p_2,k}$ where
$$p_1=\Pt{\ba{ccc}
0 & 0 & 0\\
0 & 0 & -1\\
0 & 1 & 0
\ea}\quad
p_2=\Pt{\ba{ccc}
0 & 0 & 1\\
0 & 0 & 0\\
-1 & 0 & 0
\ea}\quad
k=\Pt{\ba{ccc}
0 & -1 & 0\\
1 & 0 & 0\\
0 & 0 & 0
\ea}$$
whose commutation relations are $[p_1,p_2]=k\quad[p_2,k]=p_1\quad[k,p_1]=p_2$. Recall that $so(3)$ and $su(2)$ are isomorphic as Lie algebras, while $SU(2)$ is a double covering of $SO(3)$.

For $so(3)$ we have $\Kil(X,Y)=\Tr(XY)$ so, in particular,
$\Kil(p_i,p_j)=-2\de_{ij}$. The choice of the subspaces $$\k=\span{k}\qquad \p=\span{p_1,p_2}$$
gives a {\it Cartan decomposition} for $so(3)$.
Moreover, $\Pg{p_1,p_2}$ is an orthonormal frame for the inner product $<\cdot,\cdot>=-\frac{1}{2}\Kil(\cdot,\cdot)$ restricted to $\p$.

Defining $\Delta(g)=g\p$ and ${\mathbf g}_g(v_1,v_2)=< g^{-1}v_1,g^{-1}v_2>$, we have that $(SO(3),\Delta,{\mathbf g})$ is a $\k\oplus\p$ sub-Riemannian manifold. As for $SU(2)$, all the $\k\oplus\p$ structures that one can define on $SO(3)$ are equivalent.
\subsubsection{Expression of geodesics}
Consider an initial covector $\lam=\lam(\th,c)=\cos(\th)p_1+\sin(\th)p_2+ck\in \Lambda_\Id$: using formula \r{k+p}, we have that the exponential map is \bqn\EXP(\th,c,t)&:=&\EXP(\lam(\th,c),t)=e^{(\cos(\th)p_1+\sin(\th)p_2+ck)t} e^{-ck t}=\eqnn
$$=\Pt{\scriptsize\ba{ccc}
 K_1\cos(ct)+K_2\cos(2\th+ct)+K_3c\sin(ct)& K_1\sin(ct)+ K_2\sin(2\th+ct)-K_3 c\cos(ct)& K_4\cos(\th)+K_3\sin(\th)\\ -K_1\sin(ct)+K_2\sin(2\th+ct)+K_3c\cos(ct)& K_1\cos(ct)- K_2\cos(2\th+ct)+K_3 c\sin(ct)&-K_3\cos(\th)+K_4\sin(\th)\\
K_4\cos(\th+ct)-K_3\sin(\th+ct)&	   K_3\cos(\th+ct)+K_4\sin(\th+ct)&	\frac{\cos\Pt{\sqrt{1 + c^2} t} +c^2}{1+c^2}
\ea}$$
with $K_1=\frac{1+\Pt{1+2c^2}\cos\Pt{\sqrt{1 + c^2} t}}{2\Pt{1+c^2}}$,
$K_2=\frac{1-\cos\Pt{\sqrt{1 + c^2} t}}{2\Pt{1+c^2}}$,
$K_3=\frac{\sin\Pt{\sqrt{1 + c^2} t}}{\sqrt{1+c^2}}$,
$K_4=\frac{c\Pt{1-\cos\Pt{\sqrt{1 + c^2} t}}}{1+c^2}$.

The set of geodesics has symmetry properties similar to the $SU(2)$ case. The conjugate locus can be obtained from the one of the $SU(2)$ by the canonical projection $SU(2)\rightarrow SO(3)$. As for $SU(2)$, all the geodesics have a countable number of conjugate points.

\subsection{The $\k\oplus\p$ problem on $SL(2)$}
\label{ss-SL2}

The Lie group $SL(2)$ is the group of $2\times 2$ real matrices with determinant 1
$$SL(2)=\Pg{g\in\Mat(2,\R)\ |\ \det(g)=1}.$$
It is a non-compact group and its fundamental group is $\Z$. The Lie algebra of $SL(2)$ is the algebra of traceless $2\times 2$ real matrices $$sl(2)=\Pg{\Pt{\begin{array}{cc}
a & b\\
c & -a
\end{array}}\in\Mat(2,\R)}.$$
A basis of $sl(2)$ is $\Pg{p_1,p_2,k}$ where 
$$p_1=\frac{1}{2}\Pt{\begin{array}{cc}
1 & 0\\
0 & -1
\end{array}}\quad
p_2=\frac{1}{2}\Pt{\begin{array}{cc}
0 & 1\\
1 & 0
\end{array}}\quad
k=\frac{1}{2}\Pt{\begin{array}{cc}
0 & -1\\
1 & 0
\end{array}}$$
whose commutation relations are $[p_1,p_2]=-k\quad[p_2,k]=p_1\quad[k,p_1]=p_2$. For $sl(n)$ we have $\Kil(X,Y)=2n\Tr(XY)$, see \cite{helgason}; hence for $sl(2)$  $\Kil(X,Y)=4\Tr(XY)$ and, in particular, $\Kil(p_i,p_j)=2\de_{ij}$. The choice of the subspaces $$\k=\span{k}\qquad \p=\span{p_1,p_2}$$ provides a {\it Cartan decomposition} for $sl(2)$. For $sl(2)$ the Cartan decomposition is unique, since $\k$ must be the maximal compact subalgebra.
Moreover, $\Pg{p_1,p_2}$ is a orthonormal frame for the inner product $<\cdot,\cdot>=\frac{1}{2}\Kil(\cdot,\cdot)$ restricted to $\p$.

Defining $\Delta(g)=g\p$ and ${\mathbf g}_g(v_1,v_2)=< g^{-1}v_1,g^{-1}v_2>$, we have that $(SL(2),\Delta,{\mathbf g})$ is a $\k\oplus\p$ sub-Riemannian manifold.
\subsubsection{Expression of geodesics}
Consider an initial covector $\lam=\lam(\th,c)=\cos(\th)p_1+\sin(\th)p_2+ck\in \Lambda_\Id$: using formula \r{k+p}, we have that the exponential map is \bqn\EXP(\th,c,t)&:=&\EXP(\lam(\th,c),t)=e^{(\cos(\th)p_1+\sin(\th)p_2+ck)t} e^{-ck t}=\nn\\
&=&\scriptsize\Pt{\ba{cc}
K_1\cos\Pt{c\frac{t}{2}}+K_2\Pt{\cos\Pt{\th+c\frac{t}{2}} +c\sin\Pt{c\frac{t}{2}}}&
K_1\sin\Pt{c\frac{t}{2}}+K_2\Pt{\sin\Pt{\th+c\frac{t}{2}} -c\cos\Pt{c\frac{t}{2}}}\\
-K_1\sin\Pt{c\frac{t}{2}}+K_2\Pt{\sin\Pt{\th+c\frac{t}{2}} +c\cos\Pt{c\frac{t}{2}}}&
K_1\cos\Pt{c\frac{t}{2}}+K_2\Pt{-\cos\Pt{\th+c\frac{t}{2}} +c\sin\Pt{c\frac{t}{2}}}
\ea}
\eqnn
with
\bqn
K_1&=&\begin{graffa3}
\cosh\Pt{\sqrt{1-c^2}\frac{t}{2}}& &c\in[-1,1]\\
\cos\Pt{\sqrt{c^2-1}\frac{t}{2}}& &c\in (-\infty,-1)\cup(1,+\infty)
\end{graffa3},\nn\\
K_2&=&\begin{graffa3}
\frac{\sinh\Pt{\sqrt{1-c^2}\frac{t}{2}}}{\sqrt{1-c^2}}& &c\in(-1,1)\\
\frac{t}{2}& &c\in\Pg{-1,1}\\
\frac{\sin\Pt{\sqrt{c^2-1}\frac{t}{2}}}{\sqrt{c^2-1}}& &c\in (-\infty,-1)\cup(1,+\infty)
\end{graffa3}.
\eqnn
\subsubsection{A useful decomposition of $SL(2)$}
\label{sss:SL2dec}
\bp
For every $g\in SL(2)$ there exists a unique pair $r\in e^\k,~s\in e^\p$ such that $g=rs$.
\ep
\proof
First, notice that $e^\k=SO(2)$ and $e^\p$ is the set of $2 \times 2$ symmetric matrices with determinant 1 and positive trace.

Take $r=\Pt{\ba{cc} \cos(\th)&-\sin(\th)\\\sin(\th)&\cos(\th)\ea}\in e^\k$ and $g=\Pt{\ba{cc} \al+\de&\beta-\ga\\\beta+\ga&\al-\de\ea}\in SL(2)$. Notice that $(\al,\ga)\neq(0,0)$. We have to prove that exists a unique $\th\in\R/2\pi$ such that $s=r^{-1}g$ is symmetric with positive trace. By direct computation one gets that $s$ is symmetric if and only if $\al\sin(\th)=\ga\cos(\th)$. For any $(\al,\ga)\in\R^2\backslash(0,0)$ there exist two solutions of this equation $\th_1,\th_2\in \R/2\pi$ with $\th_2=\th_1+\pi$. Thus $\Tr\Pt{\Pt{\ba{cc} \cos(\th_1)&\sin(\th_1)\\-\sin(\th_1)&\cos(\th_1)\ea}g}=-\Tr\Pt{\Pt{\ba{cc} \cos(\th_2)&\sin(\th_2)\\-\sin(\th_2)&\cos(\th_2)\ea}g}$. Observing that a symmetric matrix with determinant 1 has nonvanishing trace, either $\th_1$ or $\th_2$ provide $\Tr(s)>0$.
\eproof

Topologically $e^\k\simeq S^1$ and $e^\p\simeq \R^2$, hence $SL(2)\simeq S^1\times\R^2$. In the following, we represent $SL(2)$ as the set $\R^2\times [0,1]$ with the identification rule $(a,b,0)\sim(a,b,1)$. See Figure \ref{fig:SL2}.
\immagine[10]{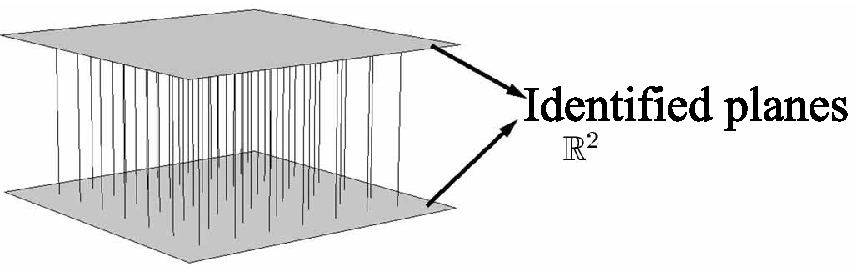}{A picture of $SL(2)$.}

\subsubsection{Symmetries in the $SL(2)$ problem}

We have the following symmetry properties:
\bi
\i cylindrical symmetry: $\Exp(\th,c,t)=e^{z_0 k}e^{x p_1 + y p_2}$ where $\Pt{\ba{c}x\\y\ea}=\Pt{\ba{cc}\cos(\th)&-\sin(\th)\\\sin(\th)&\cos(\th)\ea}\Pt{\ba{c}x_0\\y_0
\ea}$ and $(x_0,y_0,z_0)$ are defined by $\Exp(0,c,t)=e^{z_0k}e^{x_0 p_1 + y_0 p_2}$.

\i central symmetry: $\Exp(\th,-c,t)=e^{-z_0k}e^{x p_1 + y p_2}$ where $\Pt{\ba{c}x\\y\ea}=\Pt{\ba{cc}\cos(2\th)&\sin(2\th)\\\sin(2\th)&-\cos(2\th)\ea}\Pt{\ba{c}x_0\\y_0
\ea}$ and $(x_0,y_0,z_0)$ are defined by $\Exp(\th,c,t)=e^{z_0k}e^{x_0 p_1 + y_0 p_2}$.
\ei
\subsubsection{The conjugate locus}
With similar arguments to those of Section \ref{ss-SU2conj}, one checks that $g=\Exp(\th,c,t)$ is a conjugate point if and only if $$\begin{graffa3}
\sinh\Pt{\sqrt{1-c^2}\frac{t}{2}}\Pt{2\sinh\Pt{\sqrt{1-c^2}\frac{t}{2}}-t\sqrt{1-c^2}\cosh\Pt{\sqrt{1-c^2}\frac{t}{2}}}=0& &c\in(-1,1)\\
\pm \frac{t^4}{12}=0& &c=\pm 1\\
\sin\Pt{\sqrt{c^2-1}\frac{t}{2}}\Pt{2\sin\Pt{\sqrt{c^2-1}\frac{t}{2}}-t\sqrt{c^2-1}\cos\Pt{\sqrt{c^2-1}\frac{t}{2}}}=0& &c\in(-\infty,-1)\cup(1,\infty)\\
\end{graffa3}.$$

The first 2 equations have only the trivial solution $t=0$. The third one gives two series of conjugate times:
\bi
\i first series: $t_{2n-1}=\frac{2n\pi}{\sqrt{c^2-1}}$ to which correspond the conjugate loci $C^{2n-1}_\Id=e^\k\setminus\Id$;
\i second series: $t_{2n}=\frac{2x_n}{\sqrt{c^2-1}}$ where $\Pg{x_1, x_2,\ldots}$ is the ordered set of the strictly positive solutions of $x=\tan(x)$, to which correspond the conjugate loci $$C^{2n}_\Id=\Pg{\Pt{\ba{cc}
\scriptstyle \cos(x_n)\cos(y_n)+\frac{\sin\Pt{x_n}}{\sqrt{c^2-1}}\Pt{\cos\Pt{\th}+c\sin\Pt{y_n}} &
\scriptstyle \cos(x_n)\sin(y_n)+\frac{\sin\Pt{x_n}}{\sqrt{c^2-1}}\Pt{\sin\Pt{\th}-c\cos\Pt{y_n}} \\
\scriptstyle-\cos(x_n)\sin(y_n)+\frac{\sin\Pt{x_n}}{\sqrt{c^2-1}}\Pt{\sin\Pt{\th}+c\cos\Pt{y_n}} &
\scriptstyle \cos(x_n)\cos(y_n)+\frac{\sin\Pt{x_n}}{\sqrt{c^2-1}}\Pt{-\cos\Pt{\th}+c\sin\Pt{y_n}}
\ea}\mid\ba{c} c\in\R\\\th\in \R/2\pi\ea}$$
with $y_n=\frac{cx_n}{\sqrt{c^2-1}}$.
\ei
\brem
Notice that not all geodesics have conjugate points. Indeed, $\Exp(\th,c,\cdot)$ has a conjugate point if and only if $c\in(-\infty,-1)\cup(1,+\infty)$.
\erem
We present an image of the $2^{nd}$ conjugate locus (Figure \ref{fig:SL2conj}). For simplicity we present an image of its intersection with $\Pg{e^\k e^{a p_1}|a\in\R}$. The complete picture can be recovered using the cylindrical symmetry.
\immagine[10]{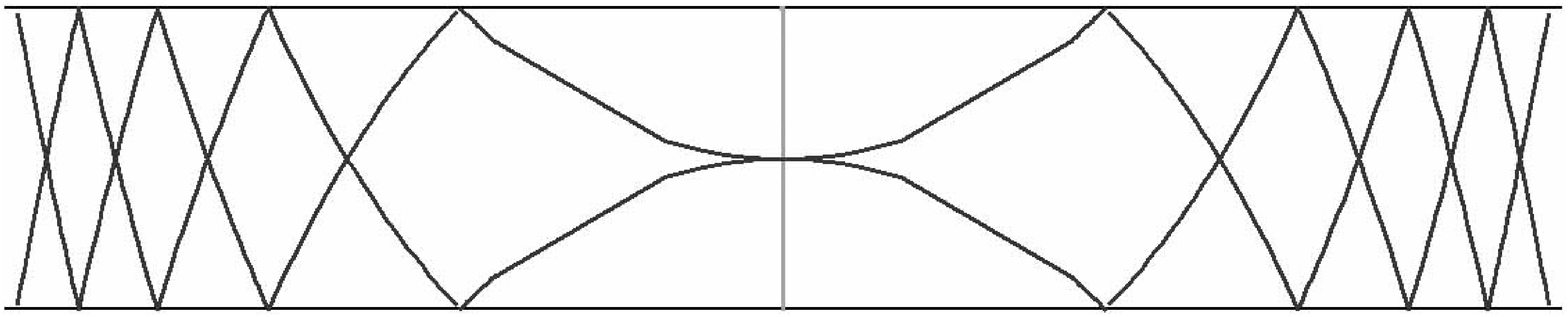}{$\k\oplus\p$ problem on $SU(2)$: section of the $2^{nd}$ conjugate locus.}
\brem
Notice that all even conjugate loci have self-intersection.
\erem
\section{A sub-Riemannian structure on lens spaces}
\llabel{s-Lpq}
\subsection{Definition of $L(p,q)$}
Fix 2 coprime integers $p,q\in\Z,~p,q\neq 0$. The \b{lens space} $L(p,q)$ is defined as the quotient of $SU(2)$ w.r.t. the identification rule $$\Pt{\ba{c}\al_1\\\beta_1\\\ea}\sim\Pt{\ba{c}\al_2\\\beta_2\\\ea}\mbox{~if~$\exists~\om\in\C$ $p$-th root of unity such that }\Pt{\ba{c}\al_2\\\beta_2\\\ea}=
\Pt{\begin{array}{cc}
\om & 0\\
0 & \om^q
\end{array}}\Pt{\ba{c}\al_1\\\beta_1\\\ea}.$$
Lens spaces are 3-dimensional compact manifolds, but excepted $L(2,1)\simeq SO(3)$, they are neither Lie groups nor homogeneous spaces of $SU(2)$. The following topological equivalences hold: $\forall~p,q,k\in\Z$, $p,q$ coprime, $p,q\neq 0$ we have $L(p,q)\simeq L(p,-q)\simeq L(-p,q)\simeq L(p, q+kp)$. Lens spaces have highly non-trivial topology, for details we refer to \cite{rolfsen}.

The following theorem permits to choose a representative of $L(p,q)$ in $SU(2)$.
\bp
\label{p-Ep}Consider the set $E_p=\Pg{\Pt{\ba{cc}\al\\\beta\ea}\in SU(2)|~\Re{\al}> 0,~\frac{\Im{\al}^2}{\sin\Pt{\frac{\pi}{p}}^2}+|\beta|^2< 1}\subset SU(2)$ and define $\partial E_p^+=\partial E_p\cap \Pg{\Im{\al}\geq 0}$, $\partial E_p^-=\partial E_p\cap \Pg{\Im{\al}\leq 0}$. Endow $\overline{E_p}$ with the equivalence relation $\div$ defined as follows:
\be
\i the relation is reflexive;
\i moreover, given $\Pt{\ba{c}\al^+\\\beta^+\ea}\in \partial E_p^+$ and $\Pt{\ba{c}\al^-\\\beta^-\ea}\in \partial E_p^-$, we have $\Pt{\ba{c}\al^+\\\beta^+\ea}\div\Pt{\ba{c}\al^-\\\beta^-\ea}$ if
\bi
\i either: $\Im{\al^+}=-\Im{\al^-}\neq 0$ and $\beta^+=e^{2\pi i\frac{q}{p}}\beta^-$;\\
\i or: $\Im{\al^+}=\Im{\al^-}=0$ and $\beta^+=e^{2\pi i\frac{n}{p}}\beta^-$ for some $n\in\Pg{1,\ldots,p}$.
\ei
\ee
The manifold $\overline{E_p}/_\div$ is diffeomorphic to $L(p,q)$.\ep
\proof
Take $\Pt{\ba{c}\al\\\beta\ea}\in SU(2)$ and let us look for $\omega$ p-th root of unity such that  $\Pt{\ba{c}\omega\al\\\omega^q\beta\ea}\in \overline{E_p}$. This condition is equivalent to
\bqn
\llabel{eq-sferaEp}
\Re{\om \al}\geq 0 &\mbox{~and~}&\frac{\Im{\om\al}^2}{\sin\Pt{\frac{\pi}{p}}^2}+|\om^q\beta|^2\leq 1.\eqn
Recalling that $|\om^q \beta|^2=|\beta|^2=1-|\al|^2$ and that $\Im{\om\al}=|\al|\sin(\arg(\om\al))$ if $\al\neq 0$, equation \r{eq-sferaEp} is equivalent to
\bqn
\arg\Pt{\om\al}\in\Pq{-\frac{\pi}{p},\frac{\pi}{p}}\mbox{~~~or~~~}\al=0 .
\eqn
Thus:
\bi
\i if $\al\neq 0$, there exist at least one solution $\om_1$ of $\arg\Pt{\om\al}\in\Pq{-\frac{\pi}{p},\frac{\pi}{p}}$. Moreover, we have 2 distinct solutions $\om_1,\om_2$ if and only if $\arg\Pt{\om_1\al}=-\frac{\pi}{p}$ and $\arg\Pt{\om_2\al}=\frac{\pi}{p}$. In this case $\Pt{\ba{c}\omega_1\al\\\omega_1^q\beta\ea}=\Pt{\ba{c}|\al|e^{-i\frac{\pi}{p}}\\\omega_1^q\beta\ea}$ and $\Pt{\ba{c}\omega_2\al\\\omega_2^q\beta\ea}=\Pt{\ba{c}|\al|e^{i\frac{\pi}{p}}\\\omega_2^q\beta\ea}$; observe that $\Pt{\ba{c}\omega_1\al\\\omega_1^q\beta\ea}\div\Pt{\ba{c}\omega_2\al\\\omega_2^q\beta\ea}$.
\i if $\al=0$, every $\om$ p-th root of unity satisfies $\Pt{\ba{c}0\\\omega^q\beta\ea}\in \overline{E_p}$; observe that for all the pairs $\om_1,\om_2$ we have $\Pt{\ba{c}0\\\omega_1^q\beta\ea}\div\Pt{\ba{c}0\\\omega_2^q\beta\ea}$.
\ei
Hence $\forall~\Pt{\ba{c}\al\\\beta\ea}\in SU(2)$ we have a unique $\Pq{\Pt{\ba{c}\omega\al\\\omega^q\beta\ea}}_\div\in \overline{E_p}/_\div$, i.e. the function  \mmfunz{\psi}{L(p,q)=SU(2)/_\sim}{\overline{E_p}/_\div}{\Pq{\Pt{\ba{c}\al\\\beta\ea}}}{\Pq{\Pt{\ba{c}\omega\al\\\omega^q\beta\ea}}_\div} is bijective.
\eproof

\brem
\label{rem:PiEp}
A crucial observation for what follows is that the projection \mmfunz{\Pi}{SU(2)}{L(p,q)}{g}{[g]} is a local diffeomorphism. Moreover, $\fz{\Pi_{|_{E_p}}}{E_p}{L(p,q)\backslash\Pq{\partial E_p}}$ is a diffeomorphism. In particular, $E_p$ contains only 1 representative for each equivalence classes of $L(p,q)$; i.e. if $g,h\in E_p$ and $[g]=[h]$, then $g=h$.
\erem
\brem
Proposition \ref{p-Ep} provides a picture of $L(p,q)$: recall that $SU(2)$ is drawn as 2 balls in $\R^3$ (see Section \ref{ss-SU2picture}). Hence $\overline{E_p}\subset SU(2)$ is drawn as a closed ellipsoid inside one of the 2 balls, via the map \mfunz{\rho}{\overline{E_p}}{\overline{B_1(0)}\subset \R^3}{\Pt{\ba{c}\al\\\beta\ea}}{(\Re{\beta},\Im{\beta},\Im{\al})}. The picture of $E_p$ is $$F_p=\Pg{(x_1,x_2,x_3)\in \overline{B_1(0)}~|~x_1^2+x_2^2+\frac{x_3^2}{\sin\Pt{\frac{\pi}{p}}^2}< 1},$$ the one of $\overline{E_p}$ is $\overline{F_p}$, see Figure \ref{fig:L41}-left.
The identification $\div$ induces the following identification on $\overline{F_p}$: given $(x_1^+,x_2^+,x_3^+)\in\partial F_p^+=\partial F_p \cap\Pg{x_3\geq 0}$ and $(x_1^-,x_2^-,x_3^-)\in\partial F_p^-=\partial F_p\cap\Pg{x_3\leq 0}$, they are identified when $x_3^+=-x_3^-$ and $\Pt{\begin{array}{c}
x_1^+\\
x_2^+
\end{array}}=
\Pt{\begin{array}{cc}
\cos(\th) & -\sin(\th)\\
\sin(\th) & \cos(\th)
\end{array}}
\Pt{\begin{array}{c}
x_1^-\\
x_2^-\\
\end{array}}\mbox{~with~}\th=\frac{2\pi q}{p}$, see Figure \ref{fig:L41}-right.
\immagine{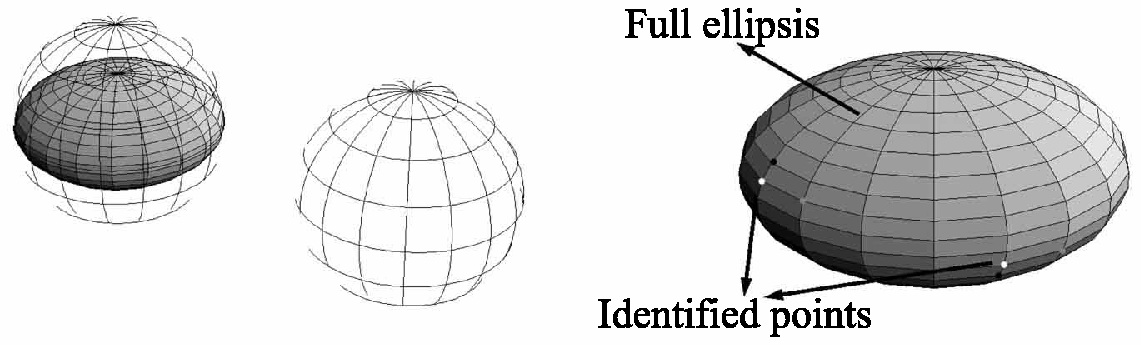}{Left: $\overline{F}_4$. Right: the representation of $L(4,1)$, with some examples of the identification rule.}
\erem
\brem
Observe that the identification rule on $\overline{F}_p$ gives a 1-to-1 identification between $\partial F_p\cap\Pg{x_3>0}$ and $\partial F_p\cap\Pg{x_3<0}$, while there are in general more identified points on $\Pg{x_1^2+x_2^2=1}\cap\Pg{x_3=0}$, see Figure \ref{fig:L41diam}.
\immagine[5]{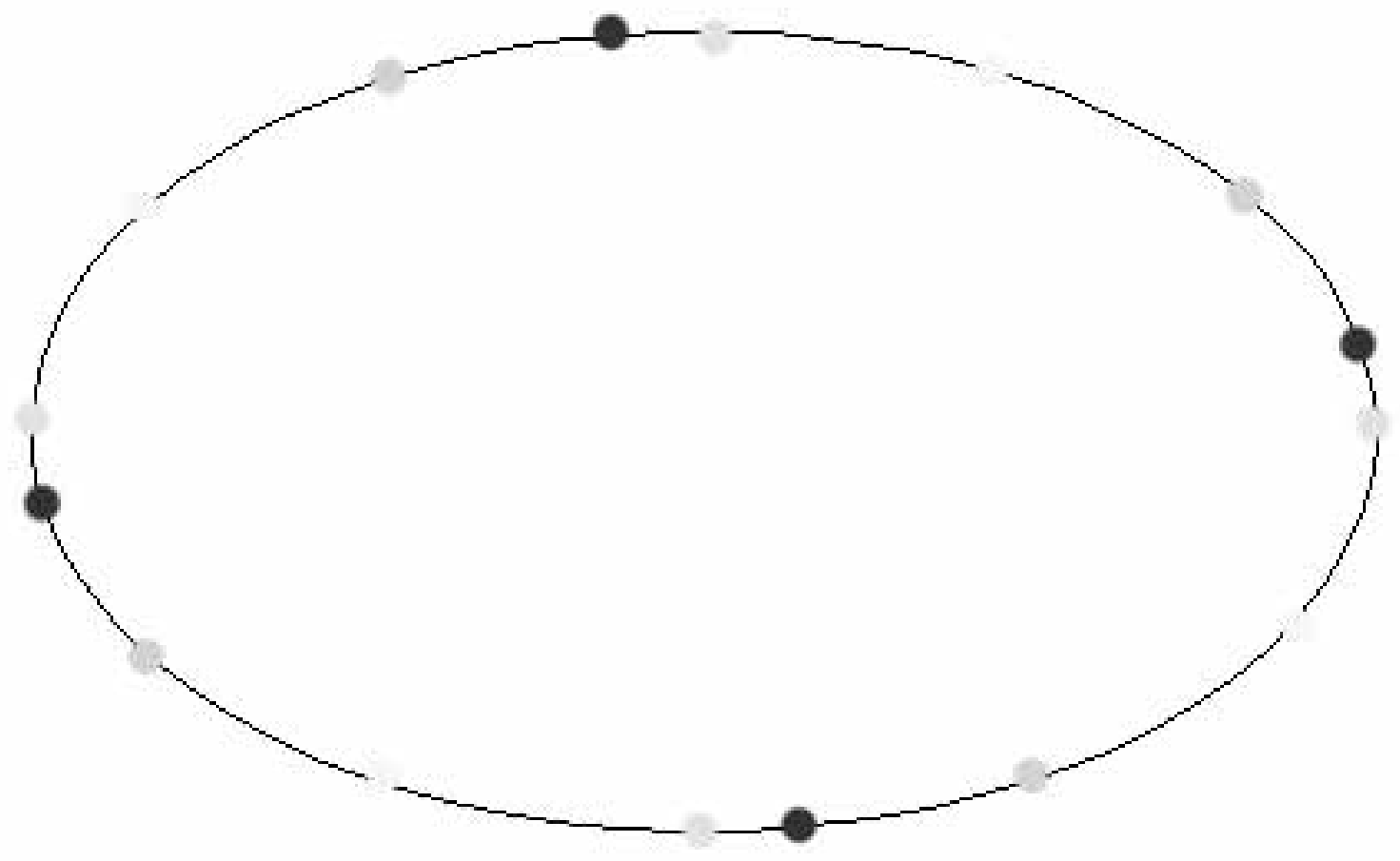}{$L(4,1)$: some examples of the identification rule on $\Pg{x_1^2+x_2^2=1}\cap\Pg{x_3=0}$.}
\erem
\subsection{Sub-Riemannian quotient structure on $L(p,q)$}
\bp
\llabel{p-projSU2Lpq}
The sub-Riemannian structure on $SU(2)$ given in Section \ref{ss-SU2} induces a 2-dim sub-Riemannian structure on $L(p,q)=SU(2)/_\sim$ via the quotient map
\mmfunz{\Pi}{SU(2)}{L(p,q)}{x}{[x]}
i.e.
\bi
\i the map $$\widetilde{\Delta}:[g]\mapsto\Pi_*\Pt{\Delta(h)}\subset T_{[g]}L(p,q)\mbox{~with~}h\in[g]$$ is a 2-dim smooth distribution on $L(p,q)$ that is Lie bracket generating;
\i $\tilde{\g}_{[g]}(v_*,w_*)=<v_*,w_*>_{[g]}:=<v,w>_h\mbox{~with~} h\in[g],~v,w\in T_hSU(2),~\Pi_*(v)=v_*,~\Pi_*(w)=w_*$ is a smooth positive definite scalar product on $\widetilde{\Delta}$.
\ei
\ep
\proof
The role of the map $\Pi$ and $\Pi_{*_{|_g}}$ is illustrated in the following diagram
$$\xymatrix{T_gSU(2) \ar[d] \ar[r]^{\Pi_{*_{|_g}}} & T_{[g]}L(p,q) \ar[d] \\
SU(2) \ar[r]^\Pi &  L(p,q).}$$

The map $\Pi$ is a local diffeomorphism, thus $\Pi_{*_{|_g}}:T_gSU(2)\to T_{[g]}L(p,q)$ is a linear isomorphism, hence $\Pi_{*_{|_g}}\Pt{\Delta(g)}$ is a 2-dim subspace of $T_{[g]}L(p,q)$.

The following statements:
\bi
\i the distribution $\widetilde{\Delta}([g])$ is well defined, i.e. $\forall h_1,h_2\in[g]$ we have $\Pi_{*_{|_{h_1}}}\Pt{\Delta(h_1)}=\Pi_{*_{|_{h_2}}}\Pt{\Delta(h_2)}$
\i the positive definite scalar product $<v_*,w_*>_{[g]}$ is well defined, i.e. $\forall h_1,h_2\in[g],~~v_1,w_1\in T_{h_1}SU(2),~~v_2,w_2\in T_{h_2}SU(2)$ such that $\Pi_{*_{|_{h_1}}}(v_1)=\Pi_{*_{|_{h_2}}}(v_2)$ and $\Pi_{*_{|_{h_1}}}(w_1)=\Pi_{*_{|_{h_2}}}(w_2)$ we have $<v_1,w_1>_{h_1}=<v_2,w_2>_{h_2}$\ei
are consequences of the following lemma:
\bl Let $h_1,h_2\in[g]$ with $h_2=\Pt{\ba{cc}\om &0\\ 0& \om^q \ea}h_1$. The map
\mmfunz{\phi}{\p}{\p}{\Pt{\ba{c} n_1\\m_1\\0\ea}}{\Pt{\ba{ccc} \Re{\omega^{q-1}}&-\Im{\omega^{q-1}}&0\\
\Im{\omega^{q-1}}&\Re{\omega^{q-1}}&0\\ 0&0&0\ea} \Pt{\ba{c} n_1\\m_1\\0\ea}}
is bijective. Moreover, it is an isometry w.r.t. the positive definite scalar product $<~,~>$ and satisfies $\forall~\eta\in\p$ $$\frac{d}{dt}_{|_{t=0}}\Pq{h_1e^{t\eta}}=\frac{d}{dt}_{|_{t=0}}\Pq{h_2e^{t\phi(\eta)}}.$$
\el
\proof
Let $h=\Pt{\ba{c}a\\b\ea}\in SU(2)$ and $\eta=(n,m,0)\in\p$. We have 
\bqn
he^{t\eta}=\Pt{
\ba{c}
a\cos\Pt{\sqrt{n^2+m^2}\frac{t}{2}}-b\sin\Pt{\sqrt{n^2+m^2}\frac{t}{2}}\frac{n-im}{\sqrt{n^2+m^2}}\\
b\cos\Pt{\sqrt{n^2+m^2}\frac{t}{2}}+a\sin\Pt{\sqrt{n^2+m^2}\frac{t}{2}}\frac{n+im}{\sqrt{n^2+m^2}}
\ea}.
\eqn

Take $h_1,h_2\in [g]$ with $h_2=\Pt{\ba{cc}\om &0\\ 0& \om^q \ea}h_1$ and $\eta_1,\eta_2\in \p$ with coordinates $\eta_1=(n_1,m_1,0)$ and $\eta_2=(n_2,m_2,0)$. Consider the trajectories $[h_1 e^{t\eta_1}]=\Pq{\Pt{
\ba{c}
a\cos\Pt{\sqrt{n_1^2+m_1^2}\frac{t}{2}}-b\sin\Pt{\sqrt{n_1^2+m_1^2}\frac{t}{2}}\frac{n_1-im_1}{\sqrt{n_1^2+m_1^2}}\\
b\cos\Pt{\sqrt{n_1^2+m_1^2}\frac{t}{2}}+a\sin\Pt{\sqrt{n_1^2+m_1^2}\frac{t}{2}}\frac{n_1+im_1}{\sqrt{n_1^2+m_1^2}}
\ea}}$ and
\bqn[h_2 e^{t\eta_2}]&=&\Pq{\Pt{
\ba{c}
\omega a\cos\Pt{\sqrt{n_2^2+m_2^2}\frac{t}{2}}-\omega^q b\sin\Pt{\sqrt{n_2^2+m_2^2}\frac{t}{2}}\frac{n_2-im_2}{\sqrt{n_2^2+m_2^2}}\\
\omega^q b\cos\Pt{\sqrt{n_2^2+m_2^2}\frac{t}{2}}+\om a\sin\Pt{\sqrt{n_2^2+m_2^2}\frac{t}{2}}\frac{n_2+im_2}{\sqrt{n_2^2+m_2^2}}
\ea}}=\nn\\
&=&\Pq{\Pt{
\ba{c}
a\cos\Pt{\sqrt{n_2^2+m_2^2}\frac{t}{2}}- \omega^{q-1} b\sin\Pt{\sqrt{n_2^2+m_2^2}\frac{t}{2}}\frac{n_2-im_2}{\sqrt{n_2^2+m_2^2}}\\
b\cos\Pt{\sqrt{n_2^2+m_2^2}\frac{t}{2}}+ \omega^{1-q} a\sin\Pt{\sqrt{n_2^2+m_2^2}\frac{t}{2}}\frac{n_2+im_2}{\sqrt{n_2^2+m_2^2}}
\ea}}.
\eqnn

Thus $\frac{d}{dt}_{|_{t=0}}[h_1e^{t\eta_1}]=\frac{d}{dt}_{|_{t=0}}[h_2 e^{t\eta_2}]$ in the case
$\begin{graffa}
n_1^2+m_1^2=n_2^2+m_2^2\\
n_1-im_1=\omega^{q-1}(n_2-im_2)\\
n_1+im_1=\omega^{1-q}(n_2+im_2)\\
\end{graffa}$, that is equivalent to $\omega^{q-1}(n_1+im_1)=n_2+im_2$. This equation is verified for $\eta_2=\phi(\eta_1)$.
\eproof\\

Since $\Pi$ is a local diffeomorphism, $\forall~g\in SU(2)~~~\exists~B(g)$ such that the map $\Pi_{*_{|_{B(g)}}}:T_{B(g)}SU(2)\to T_{B([g])}L(p,q)$ is a diffeomorphism, thus $\widetilde{\Delta}$ is smooth, Lie bracket generating, and $<v_*,w_*>_{[g]}$ is smooth as a function of $[g]$.
\eproof

Proposition \ref{p-projSU2Lpq} implies that the sub-Riemannian structures on $SU(2)$ and $L(p,q)$ defined above are locally isometric via the map $\Pi$. As a consequence, the geodesics of $(L(p,q),\tilde{\Delta},\tilde{\g})$ are the projection of geodesics of $(SU(2),\Delta,\g)$. The conjugate locus for $L(p,q)$ can be obtained from the one of $SU(2)$ by the projection $\Pi$.

\brem
One can check that the sub-Riemannian structure induced by $SU(2)$ on $L(2,1)\simeq SO(3)$ is equivalent to the $\k\oplus\p$ sub-Riemannian structure on $SO(3)$ defined in Section \ref{ss-SO3}.
\erem

\section{Cut loci and distances}
\llabel{s-cut}
In this section we prove the main theorems of the paper, i.e. we compute cut loci for $SU(2)$, $SO(3)$, lens spaces, $SL(2)$ and we prove the formula \r{eq-SU2distsystem} for the sub-Riemannian distance on $SU(2)$.

Recall that our problems satisfy the following assumptions:
{\bf i)} each point of $M$ is reached by an optimal geodesic starting from $\Id$, see Section \ref{s-k+p}; {\bf ii)} we are in the 3-D contact case, thus there are no abnormal minimizers. Hence Remark \ref{rem-agra} applies.

\bp
\llabel{t-dist}
Let $T(\th,c)$ be the cut time for $\Exp(\th,c,\cdot)$ (possibly $+\infty$ if $\Exp(\th,c,\cdot)$ is optimal on $[0,+\infty)$). Define $\D=\Pg{(\th,c,t)\in\Lam_\Id\times\R^+~|~0<t<T(\th,c)}$ and $M'=M\backslash\Pt{K_\Id\cup\Id}$. The function $\fz{\Exp_{|_\D}}{\D}{M'}$ is a diffeomorphism from $\D$ to $M'$.
\ep
\proof
Let us first check that $\Exp(\D)\subset M'$. By contradiction, let $\Exp(\th,c,t)\in M\backslash M'$, thus either $t=0$ or $t=T(\th,c)$ or $\Exp(\th,c,\cdot)$ is not optimal in $[0,t]$, i.e. $t>T(\th,c)$. Contradiction.
Let us verify that $\Exp_{|_\D}$ is injective: by contradiction, let $\Exp(\th_1,c_1,t_1)=\Exp(\th_2,c_2,t_2)$ with $(\th_1,c_1,t_1)\neq(\th_2,c_2,t_2)$. If $t_1\neq t_2$, one of the two geodesics $\Exp(\th_1,c_1,\cdot),\Exp(\th_2,c_2,\cdot)$ has already lost optimality, thus $t_i\geq T(\th_i,c_i)$, hence $\Pt{\th_i,c_i,t_i}\not\in\D$, contradiction. If $t_1=t_2$, we have that $\Exp(\th_1,c_1,t_1)$ is a cut point, hence $t_1\geq T(\th_1,c_1)$, contradiction.
To verify that $\Exp_{|_\D}$ is surjective, take $g\in M'$ and observe that there is an optimal geodesic $\Exp(\th,c,\cdot)$ reaching it at time $t\leq T(\th,c)$. But $t=T(\th,c)$ implies $g\in K_\Id$, thus $t< T(\th,c)$.

The smoothness of $\Exp_{|_\D}$ and of its inverse follows from the fact that $\Exp$ is a local diffeomorphism outside the critical points (i.e. points where the differential of $\Exp$ is not of full rank) and that the critical points do not belong to $\D$. Indeed, by contradiction, let $(\th,c,t)\in\D$ be a critical point, hence $t$ is a conjugate time: it is either the first conjugate time, that coincide with the cut time (i.e. $t=T(\th,c)$) or a greater one (i.e. $t>T(\th,c)$). In both cases $(\th,c,t)\not\in\D$. Contradiction.
\eproof

\subsection{The cut locus for $SU(2)$}
\begin{teo}
The cut locus for the $\k\oplus \p$ problem on $SU(2)$ is $$K_\Id=e^\k\setminus\Id=\Pg{e^{ck}\ |\ c\in(0,4\pi)}.$$
\end{teo}
\begin{dimostr}
Let $g\in e^\k\setminus\Id=\Pg{\Pt{\ba{c}\al\\0\ea}\ |\ \al\in\C,\ |\al|=1,\ \al\neq 1}$ and $\Exp(\th,c,\cdot)$ the minimizing geodesic steering $\Id$ to $g$ in time $T$. As a consequence of the cylindrical symmetry, we have that $\Exp(\psi,c,T)=g~~\forall~\psi\in \R/2\pi$, thus $\Pt{e^\k\backslash\Id}\subset K_\Id$.

The core of the proof is to show that there are no cut points outside $e^\k$. Recall the expression of geodesics
\bqn\EXP(\th,c,t)&=&\Pt{\ba{c}\al\\\beta\ea}\nn\\&=&\scriptsize \Pt{\begin{array}{c}
\frac{c\sin(\frac{ct}{2})\sin (\sqrt{1 + c^2}\frac{t}{2})}{\sqrt{1 + c^2}} + \cos(\frac{ct}{2})\cos(\sqrt{1 + c^2}\frac{t}{2})+i\Pt{\frac{c\cos(\frac{ct}{2})\sin (\sqrt{1 + c^2}\frac{t}{2})}{\sqrt{1 + c^2}} - \sin(\frac{ct}{2})\cos(\sqrt{1 + c^2}\frac{t}{2})}\\
\frac{\sin(\sqrt{1 + c^2}\frac{t}{2})}{\sqrt{1 + c^2}}
\Pt{\cos(\frac{ct}{2}+\th) + i\sin(\frac{ct}{2}+\th)}
\end{array}}.\eqnn
By contradiction, assume that $g\in SU(2)\backslash e^\k$ is reached by two different optimal trajectories $\Exp(\th,c,\cdot)$ and $\Exp(\psi,d,\cdot)$ at time $T$.  Observe that $\Exp(\th,c,\frac{2\pi}{\sqrt{1 + c^2}}),\Exp(\psi,d,\frac{2\pi}{\sqrt{1 + c^2}})\in e^\k\subset K_\Id$, thus
\bqn
\llabel{e-limite}
0<T<\min\Pg{\frac{2\pi}{\sqrt{1 + c^2}},\frac{2\pi}{\sqrt{1 + d^2}}}.
\eqn
Observe that $\Exp(\th,c,T)=\Exp(\psi,d,T)$ implies that $|\beta|$ is equal in the two cases, i.e.
\bqn
\llabel{eq-su2-sincsind}
\frac{\sin(\frac{\sqrt{1 + c^2} T}{2})} {\sqrt{1 + c^2}}=\frac{\sin(\frac{\sqrt{1 + d^2} T}{2})} {\sqrt{1 + c^2}}.
\eqn
From this equation it follows $|c|=|d|$. Indeed equation \r{eq-su2-sincsind} is equivalent to $\frac{\sin(\frac{\sqrt{1 + c^2} T}{2})}{\frac{\sqrt{1 + c^2}T}{2}}=\frac{\sin(\frac{\sqrt{1 + d^2} T}{2})}{\frac{\sqrt{1 + d^2}T}{2}}$. From the fact that $\frac{\sqrt{1 + c^2} T}{2},\frac{\sqrt{1 + d^2} T}{2}\in(0,\pi)$ and that the function $\frac{\sin{p}}{p}$ is injective for $p\in\Pt{0,\pi}$, it follows that $\frac{\sqrt{1 + c^2} T}{2}=\frac{\sqrt{1 + d^2} T}{2}$, hence $|c|=|d|$.

Thus we consider the two cases:
\bi
\i $c=d\in\R$: the cylindrical symmetry implies either that $\th=\psi$ (so the 2 geodesics coincide) or $g\in e^\k$. Contradiction.
\item $c=-d\in\R\backslash\Pg{0}$: with no loss of generality we assume $c>0$. Since by the central and cylindrical symmetries we have $$\Exp(\psi,-c,t)=\Pt{\ba{c}\overline{\alpha}\\e^{i(\psi+\th-\arg(\beta))}\beta\ea}
\mbox{~~~where~~~}\Exp(\th,c,t)=\Pt{\ba{c}\alpha\\\beta\ea},$$
the equation $\Exp(\th,c,t)=\Exp(\psi,-c,t)$ implies $\Im{\alpha}=0$. Hence
\bqn\llabel{e-SU2nocut} c\cos\Pt{\frac{ct}{2}}\sin\Pt{\frac{\sqrt{1+c^2}t}{2}}=
\sqrt{1+c^2}\sin\Pt{\frac{ct}{2}}\cos\Pt{\frac{\sqrt{1+c^2}t}{2}}.\eqn
The terms $c,\sin\Pt{\frac{\sqrt{1+c^2}t}{2}},\sqrt{1+c^2},\sin\Pt{\frac{ct}{2}}$ are non-zero because of equation \r{e-limite} and $c<\sqrt{1+c^2}$. Thus $\cos\Pt{\frac{ct}{2}}=0$ if and only if $\cos\Pt{\frac{\sqrt{1+c^2}t}{2}}=0$, that is impossible because $0<\frac{ct}{2}<\frac{\sqrt{1+c^2}t}{2}<\pi$.
Hence we rewrite equation \r{e-SU2nocut} as $\frac{\tan\Pt{\frac{\sqrt{1+c^2}t}{2}}}{\frac{\sqrt{1+c^2}t}{2}}=\frac{\tan\Pt{\frac{ct}{2}}}{\frac{ct}{2}}$ and state that
$\Pt{0,\tan(0)},\Pt{\frac{ct}{2},\tan\Pt{\frac{ct}{2}}},$ $\Pt{\frac{\sqrt{1+c^2}t}{2},\tan\Pt{\frac{\sqrt{1+c^2}t}{2}}}$ are 3 distinct points aligned on the graph of the function $\tan$ in $[0, \pi)$. It is impossible. Contradiction.
\ei
\end{dimostr}

The cut locus for the $\k\oplus\p$ sub-Riemannian manifold $SU(2)$ is given in Figure \ref{fig:SU2Cut}.
\immagine[9]{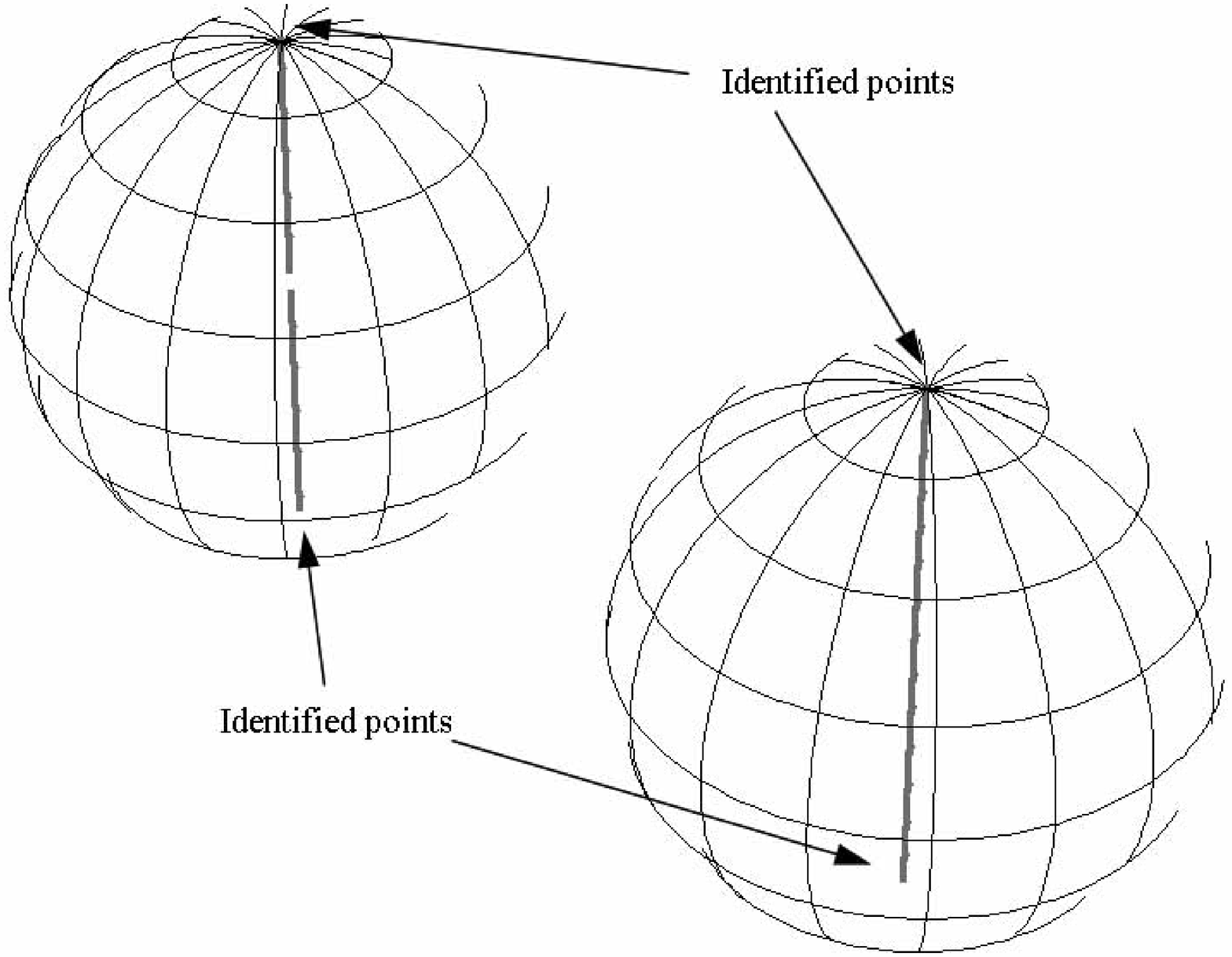}{The cut locus for the $\k\oplus\p$ sub-Riemannian manifold $SU(2)$.}
\subsubsection{The sub-Riemannian distance in $SU(2)$}
\llabel{ss-SU2dist} In this section we compute the sub-Riemannian distance on $SU(2)$, i.e. we prove Theorem \ref{t-distance}.

Let $g=\Pt{\ba{c}\al\\0\ea}=\Pt{\ba{c}e^{i\arg\Pt{\al}}\\0\ea}\in e^\k$: $g$ is reached by a geodesic $\Exp(\th,c,\cdot)$ at time $\frac{2\pi}{\sqrt{1+c^2}}$ for some $c\in\R$. Observe that $\Exp\Pt{\th,c,\frac{2\pi}{\sqrt{1+c^2}}}=\Exp\Pt{\th,\pm\sqrt{\frac{4\pi^2}{t^2}-1},t}=\Pt{\ba{c}
-\cos\Pt{\sqrt{\pi^2-\frac{t^2}{4}}}\mp i \sin\Pt{\sqrt{\pi^2-\frac{t^2}{4}}}\\ 0
\ea}=\Pt{\ba{c}
e^{i\Pt{\pi\pm\sqrt{\pi^2-\frac{t^2}{4}}}}\\ 0
\ea}$. Thus the distance $d(g,\Id)$ is the smallest $t>0$ such that $e^{i\Pt{\pi\pm\sqrt{\pi^2-\frac{t^2}{4}}}}=e^{i\arg(\al)}$, whose solution is $t=2\sqrt{\arg\Pt{\al}\Pt{2\pi-\arg\Pt{\al}}}$ where $\arg\Pt{\al}$ is chosen in $\Pq{0,2\pi}$.

Let $g=\Pt{\ba{c}\al\\\beta\ea}\in SU(2)\backslash e^\k$. Applying proposition \ref{t-dist}, we have that $\Exp_{|_\D}^{-1}(g)$ is well defined on $\D=\Pg{(\th,c,t)\in\Lam_\Id\times\R^+~|~0<t<\frac{2\pi}{\sqrt{1+c^2}}}$. Thus the sub-Riemannian distance of $g$ from the origin is $d(g,Id)=t$ where $t$ is the third component of $\Exp_{|_\D}^{-1}(g)$, i.e. the unique solution $t$ of $\Exp(\th,c,t)=g$ with $(\th,c,t)\in\D$. Using the explicit form of $\Exp$ given in \r{ss-SU2geod}, one checks that the system $\begin{graffa} \Exp(\th,c,t)=g\\ (\th,c,t)\in\D
\end{graffa}$ is equivalent to 
\bqn
\begin{graffa3}
-\frac{ct}{2}+\arctan\Pt{\frac{c}{\sqrt{1+c^2}}\tan\Pt{\frac{\sqrt{1+c^2}t}{2}}}=\arg(\al)\\
\frac{\sin\Pt{\frac{\sqrt{1+c^2}t}{2}}}{\sqrt{1+c^2}}=\sqrt{1-|\al|^2}\\
\cos(\frac{ct}{2}+\th) + i\sin(\frac{ct}{2}+\th)=\arg(\beta)
\end{graffa3}.
\eqnn
The third equation has no role for the computation of distance as a consequence of the cylindrical symmetry.

\brem The distance is a bounded function: this is due to its continuity and the compactness of $SU(2)$. The farthest point starting from $\Id$ is $-\Id$, whose distance is $2\pi$.

Notice that $\forall~\al,\beta_1,\beta_2\in\C, |\beta_1|=|\beta_2|$ we have $d\Pt{\Pt{\ba{c}\al\\ \beta_1\ea},\Id}=d\Pt{\Pt{\ba{c}\al\\ \beta_2\ea},\Id}=d\Pt{\Pt{\ba{c}\overline{\al}\\\beta_1\ea},\Id}$: this is due to the cylindrical and central symmetries.
\erem

\subsection{The cut locus for $SO(3)$ and lens spaces}

In this section we compute the cut locus for lens spaces $L(p,q)$. As a particular case, we get the cut locus for $SO(3)\simeq L(2,1)$.

\label{ss-lpcut}
\bt
The cut locus for the sub-Riemannian problem on $L(p,q)$ defined in Section \ref{p-projSU2Lpq} is a stratification
$$K_{\Pq{\Id}}=K_{\Pq{\Id}}^{sym}\cup K_{\Pq{\Id}}^{loc}$$
with
\bqn
K_{\Pq{\Id}}^{sym}&=&\Pq{\partial E_p}=\Pg{\Pq{\Pt{\ba{c} \al\\\beta\ea}}\mid a,b\in\C,\Re{\al}\geq 0,~\frac{\Im{\al}^2}{\sin\Pt{\frac{\pi}{p}}^2}+|\beta|^2= 1},\nonumber\\
K_{\Pq{\Id}}^{loc}&=&\Pq{e^\k}\setminus\Pq{\Id}=\Pg{\Pq{\Pt{\ba{c} \al\\0\ea}}\mid \al\in\C,|\al|=1,\al^p\neq 1}.
\eqnn
\et
\proof
Let us first prove the following lemma.
\bl
A geodesic $\ga(\cdot)$ in $L(p,q)$ steering $[\Id]$ to $[g]$ in minimum time $T$ admits a unique lift $\ga_0(\cdot)$ in $SU(2)$ starting from $\Id$.

Moreover, $\ga_0(t)=\Exp(\th_0,c_0,t)~~\forall~t\in[0,T]$ for some $\th_0\in\R/2\pi,~c\in\R$.
\el
\proof
Take $\ga(\cdot)$ as in the \hps. Since $L(p,q)$ and $SU(2)$ are locally diffeomorphic via $\Pi$, there is a unique lift $\ga_0(\cdot)$ in $SU(2)$ starting from $\Id$, i.e. $\ga_0(0)=\Id$ and $\Pq{\ga_0(t)}=\ga(t)~~\forall~t\in[0,T]$. 

Let us prove that $\ga_0(\cdot)$ is an optimal trajectory reaching $\ga_0(T)$. By contradiction, there exists a trajectory $\ga_1(\cdot)$ such that $\ga_1(t_1)=\ga_0(T)$ with $t_1<T$. Hence, its projection $\Pq{\ga_1(\cdot)}$ satisfies $\Pq{\ga_1(t_1)}=[g]$ with $t_1<T$. Contradiction.

Since $\ga_0(\cdot)$ is an optimal trajectory, it is a geodesic of $SU(2)$, and there exist $\th\in\R/2\pi,~c\in\R$ such that $\ga_0(t)=\Exp(\th_0,c_0,t)~~\forall~t\in[0,T]$.
\eproof\\

Let us prove that $K^{loc}_{\Pq{\Id}}\subset K_{\Pq{\Id}}$: consider $[g]\in K^{loc}_{\Pq{\Id}}$, a geodesic  steering $[\Id]$ to $[g]$ in minimum time $T$ with unique lift $\Exp(\th_0,c_0,\cdot)$. By definition of $K^{loc}_{\Pq{\Id}}$, we have $\Exp(\th_0,c_0,T)\in e^\k\backslash\Id\subset SU(2)$, i.e. $\Exp(\th_0,c_0,T)$ lies in the cut locus for the sub-Riemannian problem on $SU(2)$. Thus there exists another optimal geodesic $\Exp(\th_1,c_1,\cdot)$ defined in $[0,T]$ such that $\Exp(\th_1,c_1,T)=\Exp(\th_0,c_0,T)\in[g]$. Thus the geodesic $\Pq{\Exp(\th_1,c_1,\cdot)}$ reaches $[g]$ in minimum time. The geodesics in $SU(2)$ are distinct in a neighborhood of $\Id$, so their projections in a \nei\ of $\Pq{\Id}$ are distinct as well.

Let us now prove that $K^{sym}_{\Pq{\Id}}\subset K_{\Pq{\Id}}$: consider $[g]\in K^{sym}_{\Pq{\Id}}$, a geodesic steering $[\Id]$ to $[g]$ in minimum time $T$ with unique lift $\Exp(\th_0,c_0,\cdot)$; call $\EXP(\th_0,c_0,T)=\Pt{\ba{c}
\alpha\\
\beta\ea}\in[g]$. If $\beta=0$ we have $[g]\in K^{loc}_{\Pq{\Id}}$ or $[g]=\Pq{\Id}$, so assume $\beta\neq 0$. Due to the cylindrical and central symmetries, we have $\Exp(\th_0+\psi,-c_0,T)=\Pt{\ba{cc}
\overline{\alpha}\\
e^{2i(\th_0-\arg(\beta))+i\psi}\beta
\end{array}}$. Consider $\psi^+\in \R/2\pi$ solution of $e^{2i(\th_0-\arg(\beta))+i\psi^+}=e^{2\pi i \frac{q}{p}}$ and $\psi^-\in \R/2\pi$ solution of $e^{2i(\th_0-\arg(\beta))+i\psi^-}=e^{-2\pi i \frac{q}{p}}$. If $\EXP(\th_0,c_0,T)\in \partial E_p^+$, we have $\Pq{\Exp(\th_0+\psi^+,-c_0,T)}=\Pq{\Exp(\th_0,c_0,T)}=[g]$; if $\EXP(\th_0,c_0,T)\in \partial E_p^-$, we have similarly $\Pq{\Exp(\th_0+\psi^-,-c_0,T)}=\Pq{\Exp(\th_0,c_0,T)}=[g]$. If $c_0\neq 0$, we have found two distinct trajectories reaching $[g]$ in optimal time; if $c_0=0$, we have $\EXP(\th_0,c_0,T)\in \partial E_p^+\cap \partial E_p^-$, thus at least one among $\psi^+$ and $\psi^-$ are not null, so at least one among $\Exp(\th_0+\psi^+,0,\cdot)$ and $\Exp(\th_0+\psi^-,0,\cdot)$ are distinct from $\Exp(\th_0,0,\cdot)$ in a \nei\ of $\Id$, so are their projections in a \nei\ of $\Pq{\Id}$.

Finally, consider $[g]\in L(p,q)\backslash\Pt{K^{loc}_{\Pq{\Id}}\cup K^{sym}_{\Pq{\Id}}\cup\Pq{\Id}}$ and assume by contradiction that there exist 2 distinct geodesics steering $[\Id]$ to $[g]$ in minimum time $T$ with distinct lifts $\Exp(\th_0,c_0,\cdot),~\Exp(\th_1,c_1,\cdot)$. There are two possibilities:
\bi
\i $\Exp(\th_0,c_0,T)=\Exp(\th_1,c_1,T)$. In this case $\Exp(\th_0,c_0,T)$ lies in the cut locus for the sub-Riemannian problem on $SU(2)$, hence $[g]\in K^{loc}_{\Pq{\Id}}$. Contradiction.
\i $\Exp(\th_0,c_0,T)\neq\Exp(\th_1,c_1,T)$: since by \hp\ $[g]\not\in K^{sym}_{\Pq{\Id}}$, we have $\Exp(\th_0,c_0,T),\Exp(\th_1,c_1,T)\not\in\partial E_p$. Recalling that $[\Exp(\th_0,c_0,T)]=[\Exp(\th_0,c_0,T)]$ and $\Exp(\th_0,c_0,T),\Exp(\th_1,c_1,T)\in E_p$ implies $\Exp(\th_0,c_0,T)=\Exp(\th_1,c_1,T)$, due to remark \ref{rem:PiEp}, we have that $\Exp(\th_i,c_i,T)\in SU(2)\backslash \overline{E_p}$ for $i=0$ or $i=1$. We assume with no loss of generality that $\Exp(\th_0,c_0,T)\in SU(2)\backslash \overline{E_p}$, thus the geodesic $\Exp(\th_0,c_0,t)$ with $t\in[0,T]$ connects $\Id\in E_p$ and $\Exp(\th_0,c_0,T)\in SU(2)\backslash \overline{E_p}$; thus $\exists~\tilde{t}\in(0,T)$ such that $\Exp(\th_0,c_0,\tilde{t})\in\partial E_p$. It implies $\ga_0(\tilde{t})=\Pq{\Exp(\th_0,c_0,\tilde{t})}\in K^{sym}_{\Pq{\Id}}$, thus $\ga_0(t)$ is no more optimal for $t\in[0,T]$. Contradiction.
\ei
\eproof
\brem
Notice that $K^{loc}_{[\Id]}$ is a manifold (a circle without a point), while $K^{sym}_{[\Id]}$ is not in general. Indeed, it is an orbifold. It can be seen as $S^2\subset\R^3$ with the following identification: $(x_1^+,x_2^+,x_3^+)\in S^2\cap\Pg{x_3\geq 0}$ and $(x_1^-,x_2^-,x_3^-)\in S^2\cap\Pg{x_3\leq 0}$ are identified when $x_3^+=-x_3^-$ and $\Pt{\begin{array}{c}
x_1^+\\
x_2^+
\end{array}}=
\Pt{\begin{array}{cc}
\cos(\th) & -\sin(\th)\\
\sin(\th) & \cos(\th)
\end{array}}
\Pt{\begin{array}{c}
x_1^-\\
x_2^-\\
\end{array}}$ with $\th=\frac{2\pi q}{p}$. In the case $SO(3)\simeq L(2,1)$, we have that $K^{sym}_{[\Id]}=\R\mathbb{P}^2$ (see Figure \ref{fig:LpqCuts}-left), while in the other cases it is not locally euclidean: in fact, take a \nei\ of a point $P$ on the equator and observe that it is topologically equivalent to a set of $p$ half-planes with a common line as boundary.

Next we give an idea of the topology of the cut locus for $L(4,1)$. Consider the space $T_1$ made by the two intersecting strips $\Pg{(a,b,0)\in\R^3~|~a,b\in[-1,1]}$ and $\Pg{(a,0,b)\in\R^3~|~a,b\in[-1,1]}$ with the following identification: $(-1,b,0)\sim(1,0,b)$ and $(-1,0,b)\sim(1,-b,0)$. The boundary of this set is topologically a circle $S^1$. Consider now a 2-dimensional semisphere $T_2$. The cut locus $K^{sym}_{[\Id]}$ is topologically equivalent to the space given by glueing $T_1$ and $T_2$ along their boundaries $S^1$. The cut locus $K_{[\Id]}$ is given by glueing $K^{sym}_{[\Id]}$ with a circle $S^1$ along a point on $T_2$ and next removing a point on $S^1$ (the starting point). See a picture of it in Figure \ref{fig:LpqCuts}-right.
\erem
\immagine[12]{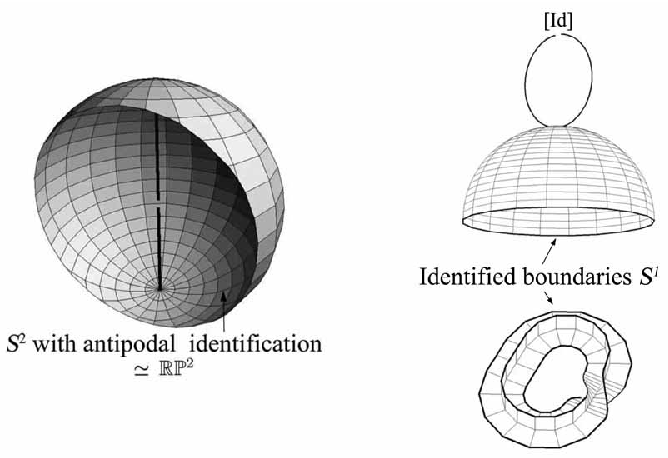}{Left: The cut locus for the sub-Riemannian problem on $SO(3)$. Right: the cut locus for the sub-Riemannian problem on $L(4,1)$.}

\subsection{The cut locus for $SL(2)$}
\bt
The cut locus for the $\k\oplus \p$ problem on $SL(2)$ is a stratification
$$K_\Id=K_\Id^{sym}\cup K_\Id^{loc}$$
with
\bqn
K_\Id^{sym}&=&e^{2\pi k}e^\p=\Pg{g\in SL(2)\mid g=g^T, \Tr{g}<0},\nonumber\\
K_\Id^{loc}&=&e^\k\setminus\Id=\Pg{\Pt{\ba{cc} \cos(\al)&-\sin(\al)\\\sin(\al)&\cos(\al)\ea}\mid \al\in R/2\pi, \al\neq 0}.
\eqnn
\et
\proof
Let us first prove that $K_\Id^{loc}\subset K_\Id$. Let $g\in e^\k\setminus\Id$: it is reached optimally by a geodesic $\Exp(\th,c,\cdot)$ at time $T$. Due to the cylindrical symmetry, we have $g=\Exp(\psi,c,T)\quad \forall\ \psi\in\R/2\pi$, thus $g\in K_\Id$.

Let us now prove that $K_\Id^{sym}\subset K_\Id$. Let $g=e^{2\pi k}e^{x_0 p_1+y_0 p_2}\in e^{2\pi k}e^\p$: it is reached optimally by a geodesic $\Exp(\th,c,\cdot)$ at time $T$.  If $x_0^2+y_0^2=0$ we have $g=e^{2\pi k}\in K_\Id^{loc}$, thus it is a cut point. If $x_0^2+y_0^2\neq 0$, due to the cylindrical and central symmetry, we have $\Exp(\th+\psi,-c,T)=e^{-2\pi k}e^{x p_1+y p_2}$ with $\Pt{\ba{c}x\\y\ea}=\Pt{\ba{cc}\cos(2\th+\psi)&\sin(2\th+\psi)\\\sin(2\th+\psi)&-\cos(2\th+\psi)\ea}\Pt{\ba{c}x_0\\y_0
\ea}$. Choose $\psi$ in such a way that $\th+\frac{\psi}{2}$ is the angle on the plane of the line passing through $(0,0)$ and $(x_0,y_0)$. In this way we have $\Pt{\ba{c}x\\y\ea}=\Pt{\ba{c}x_0\\y_0\ea}$. Observing that $e^{-2\pi k}=e^{2\pi k}$ we finally have that $g=\Exp(\th,c,T)=\Exp(\th+\psi,-c,T)$. Observe that $c\neq 0$ because $\Exp(\th,0,\cdot)\in e^\p$, thus the two geodesics $\Exp(\th,c,\cdot),\Exp(\th+\psi,-c,\cdot)$ are distinct.

We prove now that there is no cut point outside $K_\Id^{sym}\cup K_\Id^{loc}$. By contradiction, let $g\in SL(2)\setminus{\Pt{K_\Id^{sym}\cup K_\Id^{loc}\cup \Id}}$ be reached by two optimal trajectories $\Exp(\th,c,\cdot)$ and $\Exp(\psi,d,\cdot)$ at time $T$. Writing $\Exp(\th,c,t)=\Pt{\ba{cc}
g_{11}(\th,c,t)&g_{12}(\th,c,t)\\
g_{21}(\th,c,t)&g_{22}(\th,c,t)
\ea}$ we have
\bqn
\llabel{eq-r} r(c,t):=\sqrt{(g_{11}-g_{22})^2+(g_{12}+g_{21})^2}=
\begin{graffa3}
\frac{\sinh\Pt{\sqrt{1-c^2}\frac{t}{2}}}{\frac{\sqrt{1-c^2}}{2}}& &c\in\Pt{-1,1}\\
t& &c\in\Pg{-1,1}\\
\frac{\sin\Pt{\sqrt{c^2-1}\frac{t}{2}}}{\frac{\sqrt{c^2-1}}{2}}& &c\in (-\infty,-1)\cup(1,+\infty)
\end{graffa3}.
\eqn
The identity $\Exp(\th,c,T)=\Exp(\psi,d,T)$ implies $r(c,T)=r(d,T)$, that implies $c^2=d^2$. Indeed, observe that in the 3 cases described by \r{eq-r} we have respectively $r(c,t)>t,~r(c,t)=t,~r(c,t)< t$, thus either $c,d\in\Pt{-1,1}$ or $c,d\in\Pg{-1,1}$ or $c,d\in(-\infty,-1)\cup(1,+\infty)$. In each of the 3 case the identity $r(c,T)=r(d,T)$ implies $c^2=d^2$. Indeed we have the following 3 cases
\bd
\i[case $c,d\in\Pt{-1,1}$:] in this case the conclusion follows from the fact that $\frac{\sinh(p)}{p}$ is injective for $p\in\Pt{0,+\infty}$;
\i[case $c,d\in\Pg{-1,1}$:] straightforward;
\i[case $c,d\in(-\infty,-1)\cup(1,+\infty)$:] Let us prove first that $\frac{\sqrt{c^2-1} T}{2}\in(0,\pi)$. By contradiction, assume $\frac{\sqrt{c^2-1} T}{2}\geq\pi$. There exists $t\in(0,T]$ such that $\frac{\sqrt{c^2-1} t}{2}=0$, hence $r(c,t)=0$, from which it follows $\Exp(\th,c,t)\in e^\k$.

Hence either $t<T$ (and $\Exp(\th,c,\cdot)$ is not optimal on $[0,T]$, contradiction) or $t=T$ (and $g\in K^{loc}_\Id\cup\Id$, contradiction). Similarly one prove that $\frac{\sqrt{d^2-1} T}{2}\in(0,\pi)$.

Now observe that $r(c,T)=r(d,T)$ implies $\frac{\sin\Pt{\sqrt{c^2-1}\frac{T}{2}}}{\frac{\sqrt{c^2-1}T}{2}}=
\frac{\sin\Pt{\sqrt{d^2-1}\frac{T}{2}}}{\frac{\sqrt{d^2-1}T}{2}}$. 
Recalling that $\frac{\sin{p}}{p}$ is injective for $p\in(0,\pi)$ we have $\frac{\sqrt{c^2-1} T}{2}=\frac{\sqrt{d^2-1} T}{2}$, hence $c^2=d^2$.
\ed

We have two cases:
\bi
\i $c=d\in\R$. The identity $g_{11}(\th,c,T)=g_{11}(\psi,c,T)$ implies either $\th=\psi$ (i.e. the geodesics coincide) or $c\in (-\infty,-1)\cup(1,+\infty)$, thus $\sin\Pt{\sqrt{c^2-1}\frac{T}{2}}=0$, i.e. $\Exp(\th,c,T)\in e^\k$, hence either $g=\Id$ or $g\in K^{loc}_\Id$. Contradiction.
\i $c=-d\in\R\backslash\Pg{0}$. Writing $g=e^{zk}e^{x_0p_1+y_0p_2}$ we have $\Exp(\psi,-c,T)=e^{-zk}e^{xp_1+yp_2}$. The identity $\Exp(\th,c,T)=\Exp(\psi,-c,T)$ and the uniqueness of the decomposition \r{sss:SL2dec} imply $e^{zk}=\pm\Id$. Thus $g$ is symmetric, i.e. $g_{12}(\th,c,T)=g_{21}(\th,c,T)$.

If $c\in(-1,1)$ this equation implies $\frac{\tan\Pt{c\frac{T}{2}}}{c\frac{T}{2}}=\frac{\tanh\Pt{\sqrt{1-c^2}\frac{T}{2}}}{\sqrt{1-c^2}\frac{T}{2}}$. Choosing $c>0$, observe that the first positive solution $T_1$ of the equation $\frac{\tan\Pt{c\frac{T_1}{2}}}{c\frac{T_1}{2}}=\frac{\tanh\Pt{\sqrt{1-c^2}\frac{T_1}{2}}}{\sqrt{1-c^2}\frac{T_1}{2}}$ satisfies $T_1\in\Pt{\pi,3\frac{\pi}{2}}$. The other cases $c\in\Pg{-1,1}$ and $c\in (-\infty,-1)\cup(1,+\infty)$ are treated similarly and lead to $T_1\in\Pt{\pi,3\frac{\pi}{2}}$.

Thus $\cos\Pt{c\frac{T_1}{2}}<0$, hence $\Tr\Pt{g}<0$. But $\Exp(\th,c,T_1)$ symmetric and $\Tr\Pt{g}<0$ implies $\Exp(\th,c,T_1)\in K^{sym}_\Id$, i.e. $T_1$ is a cut time. Thus either $T=T_1$ (meaning that $\Exp(\th,c,t)\in K^{sym}_\Id$) or $T>T_1$ and $\Exp(\th,c,\cdot)$ is not optimal in $[0,T]$. Contradiction.
\ei
\ffoot{
\bi
\i for $c\frac{t}{2}\in\Pt{0,\frac{\pi}{2}}$ l.h.s. $>1$, while r.h.s. $<1$, thus we don't have any solution;
\i for $c\frac{t}{2}\in\Pt{\frac{\pi}{2},\pi}$ l.h.s. $<0$, while r.h.s. $>0$, thus we don't have any solution;
\i for $c\frac{t}{2}\in\Pt{\pi,3\frac{\pi}{2}}$ the l.h.s. is a continuous increasing function from 0 to $+\infty$, while the r.h.s. is continuous and increasing from $\tanh\Pt{\pi}$ to $\tanh\Pt{3\frac{\pi}{2}}$ , thus we have at least one solution of the equation: we call the first $T_1$.
\ei
For the cases $c=\pm 1$ and $c\in (-\infty,-1)\cup(1,+\infty)$ we get a different expression of the identity $g_{12}(\th,c,T)=g_{21}(\th,c,t)$ that leads nevertheless to the same result $c\frac{t}{2}\in...$}

\eproof

We give a picture of the cut locus for the $\k\oplus\p$ sub-Riemannian manifold $SL(2)$ in Figure \r{fig:SL2Cut}.
\immagine[6]{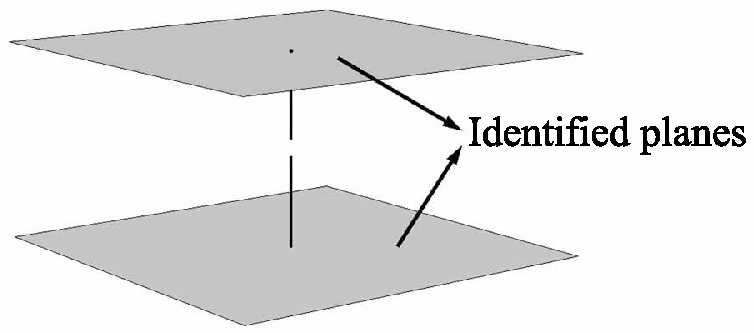}{The cut locus for the $\k\oplus\p$ sub-Riemannian manifold $SL(2)$.}

{\bf Acknowledgments:} We are grateful to A. Agrachev for many illuminating discussions. We deeply thank G. Charlot for bringing to our attention some crucial properties of the cut locus and L. Paoluzzi for many explanations on lens spaces.

\end{document}